\documentclass[a4paper]{elsarticle}
\usepackage[utf8]{inputenc}
\usepackage{import}

\usepackage{fontawesome}

\usepackage{fullpage}

\usepackage{lipsum}

\usepackage{circledsteps}
\pgfkeys{/csteps/fill color=white}
\pgfkeys{/csteps/inner ysep=5pt}
\pgfkeys{/csteps/inner xsep=5pt}

\usepackage[english]{babel}

\usepackage{enumitem}
\setlist[1]{itemsep=-5pt}
\usepackage[margin=2.45cm]{geometry}
\usepackage{amsmath} 
\usepackage{amssymb,amsthm}
\usepackage{cases}
\usepackage{mathrsfs}
\usepackage{euscript}
\usepackage{bm}

\usepackage{algorithm}
\usepackage{algpseudocode}
\theoremstyle{definition}

\usepackage{stmaryrd}
\usepackage{cases}

\usepackage{appendix}
\providecommand{\appendixname}{Appendix}
\usepackage{pdflscape}
\usepackage{tabularx}
\usepackage{booktabs}

\usepackage{multicol}
\usepackage{multirow}
\usepackage{xparse}
\usepackage{xspace}
\usepackage{makecell}
\usepackage{tablefootnote}

\usepackage{xcolor}
\usepackage{footnote} 
\usepackage{graphicx}
\usepackage{svg}
\usepackage{epstopdf}
\usepackage{tikz}
\usetikzlibrary{calc}
\usepackage{pgfplots}
\pgfplotsset{compat=newest}

\usepackage[percent]{overpic}
\usepackage{placeins} 

\usepackage[colorlinks=true]{hyperref}
\usepackage{cleveref}
\biboptions{sort&compress}
\usepackage[numbers]{natbib}
\usepackage{caption}

\usepackage[notref,notcite,final]{showkeys}

\usepackage[color=gray, backgroundcolor=yellow, textwidth=2cm, textsize=footnotesize]{todonotes}
\usepackage{setspace}
\usepackage{subcaption}
\usepackage{fullpage}
\usepackage{framed}
\usepackage[noprefix]{nomencl}

\makenomenclature

\renewcommand*\nompreamble{\begin{multicols}{2}}
\renewcommand*\nompostamble{\end{multicols}}
\usepackage{etoolbox}
\usepackage{mathtools}
\usepackage{empheq}
\usepackage{pdfpages}

\newcommand{\eps}{\bm{\varepsilon}}
\newcommand{\sig}{\bm{\sigma}}
\newcommand{\btau}{\bm{\tau}}
\newcommand{\bmu}{\bm{\mu}}
\newcommand{\gam}{\bm{\gamma}}
\newcommand{\zet}{\bm{\zeta}}

\newcommand{\Ctens}{\bm{\mathsf{C}}}
\newcommand{\Dtens}{\bm{\mathsf{D}}}

\newcommand{\Atens}{\bm{\mathsf{A}}}

\newcommand{\discrete}[1]{\mkern 1.2mu\underline{\mkern-1.2mu#1\mkern-1.2mu}\mkern 1.2mu}

\newcommand{\setD}{\mathrm{D}}
\newcommand{\setE}{\mathrm{E}}
\newcommand{\setZ}{\mathrm{Z}}

\newcommand{\Be}{\mathbf{B}^{\varepsilon}_e}
\newcommand{\Bg}{\mathbf{B}^{\gamma}_e}
\newcommand{\Bz}{\mathbf{B}^{\zeta}_e}
\newcommand{\Nc}{\mathbf{N}^{\rchi}_e}

\newcommand{\Bea}{\mathbf{B}^{\varepsilon\alpha}_e}
\newcommand{\Bga}{\mathbf{B}^{\gamma\alpha}_e}
\newcommand{\Bza}{\mathbf{B}^{\zeta\alpha}_e}
\newcommand{\Nca}{\mathbf{N}^{\rchi\alpha}_e}

\newcommand{\Bet}{\mathbf{B}^{\varepsilon\,\mathsf{T}}_e}
\newcommand{\Bgt}{\mathbf{B}^{\gamma\,\mathsf{T}}_e}
\newcommand{\Bzt}{\mathbf{B}^{\zeta\,\mathsf{T}}_e}
\newcommand{\Nct}{\mathbf{N}^{\rchi\,\mathsf{T}}_e}

\newcommand{\Beta}{\mathbf{B}^{\varepsilon\alpha\,\mathsf{T}}_e}
\newcommand{\Bgta}{\mathbf{B}^{\gamma\alpha\,\mathsf{T}}_e}
\newcommand{\Bzta}{\mathbf{B}^{\zeta\alpha\,\mathsf{T}}_e}
\newcommand{\Ncta}{\mathbf{N}^{\rchi\alpha\,\mathsf{T}}_e}

\DeclareMathOperator*{\argmin}{arg \, min}

\DeclareMathOperator*{\sym}{sym}
\DeclareMathOperator*{\dive}{div}

\newcommand{\md}{\kern0.1em---\kern0.1em\ignorespaces}
\newcommand{\nd}{\kern0.1em--\kern0.1em\ignorespaces}

\DeclareRobustCommand{\rchi}{{\mathpalette\irchi\relax}}
\newcommand{\irchi}[2]{\raisebox{\depth}{$#1\chi$}} 

\usepackage{stackengine}
\setstackgap{S}{1.75pt}

\newcommand{\trip}{\,\Shortstack{. . .}\,}

\makeatletter
\newcommand{\multiline}[1]{%
  \begin{tabularx}{\dimexpr\linewidth-\ALG@thistlm}[t]{@{}X@{}}
    #1
  \end{tabularx}
}
\makeatother

\definecolor{myblue}{RGB}{33, 33, 120}
\definecolor{mygreen}{RGB}{0, 85, 34}
\definecolor{myred}{RGB}{128, 0, 51}
\definecolor{rev}{RGB}{0, 0, 225}

\newcommand{\tB}[1]{{\color{black} #1}}

\allowdisplaybreaks

\begin{document}

\begin{frontmatter}

\author[add1]{Jacinto Ulloa\corref{cor1}}
\ead{julloa@umich.edu}
\author[add2]{Laurent Stainier}
\ead{laurent.stainier@ec-nantes.fr}

\cortext[cor1]{Corresponding author}
\address[add1]{Department of Mechanical Engineering, University of Michigan, Ann Arbor, MI 48109, USA}
\address[add2]{Nantes Universit{\'e}, Ecole Centrale Nantes, CNRS, GeM, UMR 6183, F-44000 Nantes, France}

\tnotetext[t1]{{\itshape Accepted for publication in the Journal of the Mechanics and Physics of Solids.}}

\title{
\tB{Data-driven material identification in micromorphic continua\tnoteref{t1}}
}

\begin{abstract}
We introduce a data-driven framework for identifying material behavior from full-field kinematics and \tB{external force measurements} in generalized (micromorphic) continua. \tB{The aim is to determine whether such input data can reveal generalized stress--strain states and their constitutive response without prescribing closure relations or relying on RVE-based homogenization. To this end, the approach infers the associated generalized stresses from full-field boundary value problems and constructs representative material datasets via clustering in a non-classical phase space}. We show that the proposed method reliably extracts non-symmetric and higher-order local stress states, providing material data suitable for either model calibration or model-free data-driven simulations of generalized continua. These capabilities are demonstrated in linear and nonlinear validation simulations with synthetic data, and in an application to mechanical metamaterials, suggesting a practical route for material characterization of microstructured~solids.
\end{abstract}


\end{frontmatter}

\section{Introduction}
Generalized continuum theories emerged from the need to describe material systems whose mechanical response is strongly influenced by internal microstructural effects, for which the physics supplied by classical Cauchy continuum theories is no longer sufficient~\cite{maugin2010generalized}. Starting with the seminal work of the Cosserat brothers~\cite{cosserat1909}, which gave rise to the notion of micropolar continua, several decades of research have produced a broad family of generalized continuum formulations, including microstrain~\cite{forest2006}, gradient-enhanced~\cite{triantafyllidis1986gradient}, non-local~\cite{lazar2006theory}, and micromorphic theories~\cite{eringen1968mechanics}. These theories capture complex phenomena related to microstructure~\cite{ariza2024homogenization, glaesener2019continuum, xu2024derivation, sarhil2024computational}, strain localization~\cite{triantafyllidis1986gradient, muhlhaus1991, forest2004, collins2020cosserat}, size effects~\cite{chen1998fracture, dillard2006, rubin2023eulerian,Ulloa2024b}, and wave dispersion~\cite{madeo2015, misra2016,dayal2017, shaat2018reduced} in a range of materials. Yet, despite substantial theoretical progress, the practical identification of constitutive relations for generalized continua remains a persistent challenge, particularly for complex micromorphic models, often rendering them elusive in applications. This limitation stems from the fact that generalized material states, in particular stress quantities, are inherently unmeasurable. Hence, the objective of the present work is to develop a model-free framework for identifying material stress--strain data directly from full-field kinematics and force measurements in generalized continua.

The cornerstone of Cosserat continua is the introduction of micro-rotational degrees of freedom and associated balance laws, giving rise to a continuum theory in which the classical symmetry of the stress tensor is relaxed. Non-symmetric stresses of this type arise from body and surface couples of mechanical origin or from other physical mechanisms, e.g., through the action of electromagnetic fields~\cite{ivanova2021new}. Although conceptually significant, the notion of generalized continua thus initiated remained dormant for several decades before being revived and substantially generalized in the 1960s and 1970s by~\citet{eringen1964,eringen2012microcontinuum}, \citet{mindlin1964}, and~\citet{germain1973}. These developments extended the Cosserat framework from rigid micro-rotations to fully deformable microstructures via generally non-symmetric micro-deformation tensors and third-order stress measures, thereby establishing the class of micromorphic continua. The resulting theory is notably general, embedding several models as special cases (see~\citet{neff2014} for a unified overview), including micropolar and second-gradient elasticity, and has since been adapted to describe, e.g., plasticity~\cite{forest2003,regueiro2009,bryant2019,rys2020,lindroos2022} and damage~\cite{aslan2011,miehe2017,brepols2017,yin2022} (see~\citet{forest2009} for an overview in this context). Despite its unifying structure, practical applications of phenomenological micromorphic models \tB{face challenges in experimental identification and constitutive calibration due to the large number of material parameters} relating generalized strains to generalized stresses.

Pivoting away from traditional parameter-based constitutive models, model-free data-driven methods~\cite{kirchdoerfer2016, eggersmann2019, karapiperis2021a, ulloa2023, prume2025direct} circumvent the need for an explicit constitutive parameterization of material behavior. Rather than postulating phenomenological laws, these methods close the fundamental kinematic and balance equations with empirical data obtained directly from experiments or lower-scale simulations. While this approach is mostly developed for classical Cauchy continua, recent efforts have considered extensions to generalized continuum theories~\cite{karapiperis2021b, kamasamudram2023, ulloa2024}. Nevertheless, much like constitutive models, generalized data-driven frameworks critically depend on the availability of generalized stress--strain data in the appropriate phase space, albeit intended for direct use in model-free simulations rather than parameter fitting. 

Ideally, homogenization methods, based on either formal two-scale asymptotics or discrete-to-continuum variational analysis, would provide (i) an unambiguous effective continuum (energy) form and (ii) explicit effective material parameters. Such approaches are well established for linear Cauchy continua, and have also been developed for generalized continua arising from particular microstructures~\cite{ariza2024homogenization, Ulloa2024b, glaesener2019continuum, misra2015, xu2024derivation}. However, more often, computational homogenization is employed, e.g., to numerically extract generalized continuum quantities from heterogeneous but classical (Cauchy) materials defined at the scale of representative volume elements (RVEs).  Developments in this context include homogenization toward micropolar~\cite{forest1998cosserat, feyel2003multilevel, de2011cosserat}, second-gradient~\cite{kouznetsova2002multi, kouznetsova2004multi, nguyen2013multiscale}, and micromorphic continua~\cite{forest2002homogenization, janicke2012minimal, alavi2021construction, biswas2017micromorphic, ehlers2020, rokovs2019, neff2020identification, sarhil2024computational}. More recently, \citet{miller2022micromorphic} proposed a filter-based method that extracts micromorphic stresses and deformation measures via averaging rules in fully resolved direct numerical simulations. While these approaches demonstrate that generalized stress--strain data can, in principle, be extracted from microstructured materials, most RVE-based homogenization schemes still depend critically on assumptions regarding the kinematic ansatz and higher-order boundary conditions at the unit-cell level. In turn, filter-based homogenization relies on microscale stresses and associated numerical averaging procedures. Consequently, no clear consensus currently exists on the computational characterization of effective generalized material states. It is worth noting that recent progress in this context has been made using dynamical systems theory~\cite{roberts2025accurate}.

\tB{Here, we propose to pivot away from homogenization and instead address the question of whether generalized stress--strain data can be inferred directly from full-field kinematic measurements and applied forces in boundary value problems (BVPs)}. Building upon previous developments for classical Cauchy materials~\cite{leygue2018data, stainier2019model}, the proposed framework is agnostic to the constitutive behavior and relies solely on enforcing non-classical balance laws and compatibility relations. In this way, it enables a direct identification of material data for \tB{micromorphic continua}, including non-symmetric and higher-order stresses. The identified data is then suitable for either material model calibration or model-free data-driven simulations.

The paper is structured as follows. \tB{Section~\ref{sec:micromorphic} provides an overview of micromorphic mechanics and its forward data-driven formulation}. Section~\ref{sec:DDI} presents the proposed data-driven identification method. Section~\ref{sec:num} presents synthetic validation examples, including elastic and inelastic materials, followed by identification and prediction in mechanical metamaterials. Section~\ref{sec:conclusion} provides a summary and~outlook.

\tB{
\section{Data-driven computational mechanics for generalized continua}
\label{sec:micromorphic} 

\tB{For the reader's convenience, we first provide an overview of micromorphic mechanics at small strains and the data-driven formulation of the forward problem, following the presentation in~\citet{ulloa2024}.}

\subsection{Overview of micromorphic continua at small strains}

\subsubsection{Kinematics and balance laws}
}

\begin{figure}[b!]
    \vspace{1em}
    \centering
    \small
    \includeinkscape[scale=0.7]{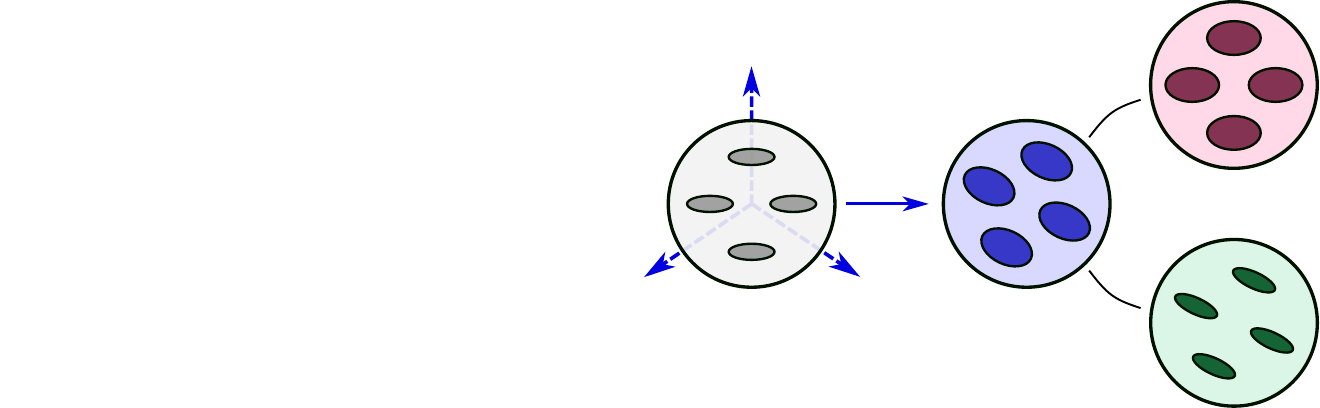}
    \vspace{1em}
    \caption{\tB{Schematic representation of a BVP in a micromorphic continuum (left) with a generic microstructure (right). The dashed blue arrows on the left represent second-order tensors. The first-order theory adopted in this work assumes that the micro-deformation tensor ${\color{blue}\bm{\rchi}} = {\color{myred}\mathrm{sym}\,\bm{\rchi}} + {\color{mygreen}\mathrm{skw}\,\bm{\rchi}}$ is homogeneous in $\Omega^\mathrm{m}$ but varies in $\Omega$.}}
    \label{fig:micbvp}
\end{figure}

Consider a \emph{microstructured} deformable solid $\Omega\subset\mathbb{R}^n$ of spatial dimension $n$. The solid is described as a \tB{\emph{micromorphic continuum}}~\cite{eringen1964,mindlin1964}, such that, attached to each macro-coordinate $\bm{x}\in\Omega$, is a \emph{micro-continuum} body $\Omega^\mathrm{m}$ defined through micro-coordinates $\bm\xi\in\Omega^\mathrm{m}$ and subject to independent deformations $\bm\rchi$ (Figure~\ref{fig:micbvp}). Then, the motion of a point $(\bm{x},\bm{\xi})$ is given by a \emph{micro-displacement} field expressed as the sum of macro-displacement and fluctuation~terms:
\begin{equation}
    \bm{v}(\bm{x},\bm{\xi})=\bm{u}(\bm{x})+\bm{\rchi}(\bm{x})\cdot\bm{\xi} + O(\Vert\bm{\xi}\Vert^2), \qquad \rchi_{ij}=\frac{\partial v_{i}}{\partial\xi_j}.
\end{equation}
The deformation process is thus characterized by the macro-displacement field $\bm{u}:\Omega\to\mathbb{R}^n$ and the micro-deformation field $\bm{\rchi}:\Omega\to\mathbb{R}^{n\times n}$, assumed homogeneous in $\Omega^\mathrm{m}$ for a first-order theory~\cite{germain1973}.

The macroscopic BVP is then given in terms of the primary fields $\bm{u}$ and $\bm{\rchi}$, the latter representing an additional set of degrees of freedom (DOFs). External actions are applied on the solid boundary $\Gamma$, decomposed into a Dirichlet part $\Gamma^u_\mathrm{D}$ with imposed displacements $\bar{\bm{u}}$ and a Neumann part $\Gamma^u_\mathrm{N}$ with applied tractions $\bar{\bm{t}}(\bm{x})\in\mathbb{R}^{n}$, where $\Gamma^u_\mathrm{D}\cup\Gamma^u_\mathrm{N}=\Gamma$ and $\Gamma^u_\mathrm{D}\cap\Gamma^u_\mathrm{N}=\emptyset$. Likewise, we may consider imposed micro-deformations $\bar{\bm{\rchi}}$ on $\Gamma^\rchi_\mathrm{D}$ and applied double tractions $\bar{\bm{T}}(\bm{x})\in\mathbb{R}^{n\times n}$ on $\Gamma^\rchi_\mathrm{N}$, with $\Gamma^\rchi_\mathrm{D}\cup\Gamma^\rchi_\mathrm{N}=\Gamma$ and $\Gamma^\rchi_\mathrm{D}\cap\Gamma^\rchi_\mathrm{N}=\emptyset$. The body may also be subjected to volume forces $\bm{b}(\bm{x})\in\mathbb{R}^{n}$ and double forces~$\bm{M}(\bm{x})\in\mathbb{R}^{n\times n}$.

The first-order micromorphic continuum involves the following generalized strain measures:
\begin{align}
&\text{Compatible strain} &\eps=\sym\,\nabla\bm{u},&&& &\varepsilon_{ij}=\frac{1}{2}(u_{i,j}+u_{j,i});&&&
\label{eq:epsilon}   \\
&\text{Relative strain}   &\gam=\nabla\bm{u} - \bm{\rchi},&&& &\gamma_{ij}=u_{i,j}-\rchi_{ij};&&&
\label{eq:gamma} \\ 
&\text{Micro-deformation gradient} &\zet=\nabla\bm{\rchi},&&& &\zeta_{ijk}=\rchi_{ij,k}.&&&
\label{eq:zeta}
\end{align}
The mechanical power density then reads
\begin{equation}
p^\mathrm{int}=\sig:\dot{\eps}+\btau:\dot{\gam}+\bmu\trip\dot{\zet}, 
\label{eq:pdens}
\end{equation}
where $\sig$ is the second-order, symmetric Cauchy stress tensor; $\btau$ is the second-order, generally non-symmetric relative stress tensor; and $\bmu$ is the third-order double stress tensor. An appeal to the principle of virtual power~\cite{germain1973} yields the following balance equations:
\begin{align}
&\text{Stress equilibrium} \qquad &\dive\,(\sig+\btau) + \bm{b} = \bm{0}\ \ \ \text{in } \Omega,&&& \quad &(\sig+\btau)\cdot\bm{n}=\bar{\bm{t}}  \ \ \ \text{on } \Gamma_{\mathrm{N}}^{u};&&&
\label{eq:eqsig}   \\
&\text{Double stress equilibrium} \qquad  &\dive\,\bmu+\btau + \bm{M} = \bm{0} \ \ \ \text{in }  \Omega,&&& \quad &\bmu\cdot\bm{n}=\bar{\bm{T}}  \ \ \ \text{on } \Gamma_{\mathrm{N}}^{\rchi}.&&&
\label{eq:eqmu}
\end{align}
This system must be solved for the primary fields, $\bm{u}\in\mathrm{H}^1(\Omega;\mathbb{R}^n)$ and $\bm{\rchi}\in\mathrm{H}^1(\Omega;\mathbb{R}^{n\times n})$, subject to boundary conditions $\bm{u}=\bar{\bm{u}}$ on $\Gamma^u_\mathrm{D}$ and $\bm{\rchi}=\bar{\bm{\rchi}}$ on $\Gamma^\rchi_\mathrm{D}$. 

\smallskip

\subsubsection{Constitutive behavior}

At this point, employing a traditional modeling approach requires closure relations between the generalized stresses and the conjugate generalized strains, i.e., a mapping $(\eps,\gam,\zet)\mapsto \{\sig,\btau,\bmu\}(\eps,\gam,\zet)$, a far-from-trivial requirement given the high-dimensional phase space. Indeed, even for isotropic elasticity, first-order micromorphic theories may require up to 18 parameters\md of course, this general model is often simplified by means of kinematic assumptions, yielding reduced theories such as the micropolar model (with $\bm\gamma$ purely skew-symmetric)~\cite{cosserat1909}, the microstrain model (with $\bm\gamma$ purely symmetric)~\cite{forest2006}, or the relaxed micromorphic model (replacing $\nabla\bm{\rchi}$ with the curvature part $\mathrm{curl}\,\bm{\rchi}$)~\cite{neff2014}. Still, as discussed in the data-driven literature for both Cauchy~\cite{kirchdoerfer2016,stainier2019model} and generalized~\cite{karapiperis2021b,ulloa2024} continua, material behavior is, in principle, only known through an empirical dataset, be it experimental or from high-fidelity lower-scale simulations. Such data may not be accurately represented by a closed-form model with a predefined functional form and material parameters, particularly for complex material behavior. 

To address this problem in the context of generalized continua, we have recently presented a data-driven approach to micromorphic mechanics~\cite{ulloa2024}, where we assume that the material response is encoded in an empirical dataset in terms of generalized stress--strain coordinates $\{(\eps,\gam,\zet),(\sig,\btau,\bmu)\}$. Such data is employed directly to solve the fundamental equations~\eqref{eq:epsilon}--\eqref{eq:zeta} and~\eqref{eq:eqsig}--\eqref{eq:eqmu} in BVPs without resorting to a constitutive model. However, the study~\cite{ulloa2024} considered synthetic datasets for proof-of-concept simulations, leaving the open question of whether the required generalized stress--strain data can be obtained from microstructured material responses in a practical scenario. The present work aims to fill this gap through a data-driven procedure for generalized material data identification, as described in section~\ref{sec:DDI}.

\subsection{\tB{Data-driven formulation}}
\label{sec:fem}

\tB{To set the stage for the identification framework, let us briefly recall the forward data-driven problem for micromorphic mechanics presented in~\citet{ulloa2024}}. Consider a discretization of $\Omega$ into $N$ nodes and $M$ material points, achieved via finite elements or another means. We assemble nodal displacements and external forces in vectors $\discrete{\mathbf{u}}=\{{\bm{u}}_a\in\mathbb{R}^{n}\}_{a=1}^N$ and $\discrete{\mathbf{f}}=\{{\bm{f}}_a\in\mathbb{R}^{n}\}_{a=1}^N$, respectively. Similarly, the nodal micro-deformations and double forces are collected in vectors $\discrete{\bm{\rchi}}=\{{\bm{\rchi}}_a\in\mathbb{R}^{n_\chi}\}_{a=1}^N$ and $\discrete{\mathbf{m}}=\{{\bm{m}}_a\in\mathbb{R}^{n_\chi}\}_{a=1}^N$, respectively. Here, the size of the nodal micro-deformation vector, $n_\chi$, depends on the symmetry conditions, e.g., $n_\chi=n^2$ for the full micromorphic model and $n_\chi=n(n-1)/2$ for the micropolar~case. 

Hereafter, second- and higher-order tensors evaluated at a point in space are expressed in Voigt form.  The discrete form of the kinematic relations~\eqref{eq:epsilon}--\eqref{eq:zeta} then reads
\begin{align}
&\eps_e = \Be\,\discrete{\mathbf{u}},
\label{eq:epsilon_disc}   \\
&\gam_e = \Bg\,\discrete{\mathbf{u}}\;-\;\Nc\,\discrete{\bm{\rchi}},
\label{eq:gamma_disc} \\ 
&\zet_e = \Bz\,\discrete{\bm{\rchi}},
\label{eq:zeta_disc}
\end{align}
while the balance equations~\eqref{eq:eqsig}--\eqref{eq:eqmu} take the discrete forms
\begin{align}
&\sum_{e=1}^{M}w_e\Big[\Bet\,\sig_e\;+\;\Bgt\,\btau_e\Big]-\discrete{\mathbf{f}}=\bm{0},
\label{eq:eqsig_disc}\\
&\sum_{e=1}^{M}w_e\Big[\Bzt\,\bmu_e\;-\;\Nct\,\btau_e\Big]-\discrete{\mathbf{m}}=\bm{0}.
\label{eq:eqmu_disc}
\end{align}
Here, $w_e\in\mathbb{R}$ is a standard integration weight; $\mathbf{N}^\rchi_e\in\mathbb{R}^{n_\chi\times N n_\chi}$ is a shape function matrix for $\discrete{\bm{\rchi}}$; and $\Be\in\mathbb{R}^{n_\varepsilon \times N n }$, $\Bg\in\mathbb{R}^{n_\chi\times N n_\chi}$, and $\Bz\in\mathbb{R}^{n_\zeta\times N n_\chi}$ are standard discretized gradient operators, with $n_\varepsilon=n(n+1)/{2}$ and $n_\zeta=n_\chi(n_\chi+1)/2$.

The state of the solid at a material point $e$ is characterized by generalized coordinates $\mathbf{z}_e$ in a local phase space $\setZ_e$. In the micromorphic setting considered here, we have
\begin{equation}
    \mathbf{z}_e\coloneqq(\eps_e,\gam_e,\zet_e,\sig_e,\btau_e,\bmu_e)\in\mathrm{Z}_e,\quad \mathrm{Z}_e\coloneqq\mathbb{R}^{n_\varepsilon}\times\mathbb{R}^{n_\chi}\times\mathbb{R}^{n_\zeta}\times\mathbb{R}^{n_\varepsilon}\times\mathbb{R}^{n_\chi}\times\mathbb{R}^{n_\zeta}.  
    \label{eq:ps}
\end{equation}
This local phase space is equipped with the metric
\begin{equation}
\Vert \mathbf{z}_e \Vert_e^2 \coloneqq \frac{1}{2}\Big(\Ctens_e\eps_e\cdot\eps_e + \Ctens_e^{-1}\sig_e\cdot\sig_e + \Dtens_e\gam_e\cdot\gam_e + \Dtens_e^{-1}\btau_e\cdot\btau_e + \Atens_e\zet_e\cdot\zet_e + \Atens_e^{-1}\bmu_e\cdot\bmu_e\Big),
\label{eq:locmet}
\end{equation}
where $\Ctens_e\in\mathbb{R}^{n_\varepsilon\times n_\varepsilon}_{\mathrm{sym},+}$, $\Dtens_e\in\mathbb{R}^{n_\chi\times n_\chi}_{\mathrm{sym},+}$, and $\Atens_e\in\mathbb{R}^{n_\zeta\times n_\zeta}_{\mathrm{sym},+}$ are positive-definite matrices representing numerical metric tensors (\ref{sec:metric}) in Voigt form. The local coordinates, collected as $\mathbf{z}\coloneqq\{\mathbf{z}_e\in\mathrm{Z}_e\}_{e=1}^M\in\mathrm{Z}$, represent points in the global phase space $\mathrm{Z}=\mathrm{Z}_1\times\dots\times\mathrm{Z}_M$. The global distance metric then reads
\begin{equation}
    \Vert \mathbf{z} \Vert^2 \coloneqq \sum_{e=1}^M w_e \Vert \mathbf{z}_e \Vert_e^2.
    \label{eq:globmet}
\end{equation}

With these definitions in hand, the data-driven micromorphic problem is given by
\begin{equation}
        \min_{\displaystyle\bar{\mathbf{z}}\in\setD}\,\min_{\displaystyle\mathbf{z}\in\setE} \Vert \bar{\mathbf{z}}-\mathbf{z}\Vert^2 = \min_{\displaystyle\mathbf{z}\in\setE}\,\min_{\displaystyle\bar{\mathbf{z}}\in\setD} \Vert \bar{\mathbf{z}}-\mathbf{z}\Vert^2 ,
\label{eq:DDproblem}
\end{equation}
where $\setE\subset\mathrm{Z}$ is the \emph{generalized equilibrium set}
\begin{equation}
\setE\coloneqq\big\{\mathbf{z}\in\mathrm{Z} \ \colon \  \text{\eqref{eq:epsilon_disc}--\eqref{eq:eqmu_disc}}\big\}.
\label{eq:Eset}
\end{equation}
On the other hand, $\setD=\setD_1\times\dots\times\setD_M$ is a \tB{\emph{given}} global dataset, where $\setD_e\subset\mathrm{Z}_e$ denotes the local dataset containing generalized stress--strain states accessible to the material point $e$.

\tB{
\section{Data-driven identification for generalized continua}
\label{sec:DDI}

We now turn from the forward data-driven problem recalled above to the inverse problem of data-driven identification considered in this work. In the forward problem~\eqref{eq:DDproblem}, as well as in model-based formulations requiring calibration, a material dataset $\setD$ in terms of $(\eps,\gam,\zet,\sig,\btau,\bmu)$ is assumed to be available. Here, instead, the available information consists of kinematic field measurements and applied forces collected from experimental or numerical BVPs with arbitrary geometries and loading conditions, from which we seek to identify the unknown generalized stress fields and material dataset $\setD$.

In this sense, the procedure provides an alternative to RVE-based homogenization for obtaining generalized stress--strain states in microstructured materials, a task for which no consensus seems to exist. The resulting data may then be used for model-free data-driven simulations in the sense of~\eqref{eq:DDproblem}, constitutive model calibration, or material characterization more broadly.}

As in data-driven identification (DDI) for Cauchy continua~\cite{leygue2018data,stainier2019model}, we suppose that full-field displacement measurements are given for a BVP under various loading conditions. It is reasonable to assume the availability of such measurements from micromechanical simulations or experimental images. In the present case of microstructured media, we further assume that the material has an identifiable representative volume $\Omega_\mathrm{m}$ over which, considering a separation of scales, micro-deformations are also available. Examples of how such kinematic measures may arise are found in the literature~\cite{forest2011, hutter2019theory} and in our example in section~\ref{sec:MMM}. Then, we may construct a database over $N_\mathrm{L}$ loading conditions indexed by $\alpha\in\{1,\dots, N_\mathrm{L}\}$, consisting of discretized displacements $\discrete{\mathbf{u}}^\alpha=\{{\bm{u}}^\alpha_a\in\mathbb{R}^{n}\}_{a=1}^N$ and micro-deformations $\discrete{\bm{\rchi}}^\alpha=\{{\bm{\rchi}}_a^\alpha\in\mathbb{R}^{n_\chi}\}_{a=1}^N$, together with applied forces $\discrete{\mathbf{f}}^\alpha=\{{\bm{f}}^\alpha_a\in\mathbb{R}^{n}\}_{a=1}^N$~and~{$\discrete{\mathbf{m}}^\alpha=\{\bm{m}_a^\alpha\in\mathbb{R}^{n_\chi}\}_{a=1}^N$}.

The objective of the proposed DDI framework is, then, to find \tB{a prescribed number} $\bar{N}$ of material data points $\bar{\mathbf{z}}_i\coloneqq(\bar{\eps}_i,\bar{\gam}_i,\bar{\zet}_i,\bar{\sig}_i,\bar{\btau}_i,\bar{\bmu}_i)$, such that:

\begin{enumerate}
    \item for each loading case $\alpha$ and material point $e$, it is possible to find generalized stresses $(\sig_e^\alpha,\btau_e^\alpha,\bmu_e^\alpha)$ satisfying the balance equations~\eqref{eq:eqsig_disc}--\eqref{eq:eqmu_disc};
    \item for each loading case $\alpha$, the generalized material state $\bar{\mathbf{z}}_{i_e^\alpha}=(\bar{\eps}_{i_e^\alpha},\bar{\gam}_{i_e^\alpha},\bar{\zet}_{i_e^\alpha},\bar{\sig}_{i_e^\alpha},\bar{\btau}_{i_e^\alpha},\bar{\bmu}_{i_e^\alpha})$ assigned to each material point $e$ minimizes the distance to the mechanically admissible state $\mathbf{z}^\alpha_e=(\eps^\alpha_e,\gam^\alpha_e,\zet^\alpha_e,\sig^\alpha_e,\btau^\alpha_e,\bmu^\alpha_e)$, according to the local metric~\eqref{eq:locmet}.
\end{enumerate}

\noindent Note that the mechanical strains $(\eps^\alpha_e,\gam^\alpha_e,\zet^\alpha_e)$ are known at every material point $e$ from the available kinematic measurements $\discrete{\mathbf{u}}^\alpha$ and $\discrete{\bm{\rchi}}^\alpha$ and the discrete kinematic relations~\eqref{eq:epsilon_disc}--\eqref{eq:zeta_disc}.

We seek the two objectives above through the solution of the minimization problem
    \begin{equation}    \min_{\displaystyle\{(\sig^\alpha_e,\btau^\alpha_e,\bmu_e^\alpha)\}}\;\min_{\displaystyle\vphantom{(}\{\bar{\mathbf{z}}_i\}}\;\min_{\displaystyle\vphantom{(}\{i_e^\alpha\}} \bigg\{\sum_{\alpha=1}^{N_\mathrm{L}}\sum_{e=1}^M w_e \Vert \mathbf{z}^\alpha_e - \bar{\mathbf{z}}_{i_e^\alpha}\Vert_e^2 \ \colon \  \text{\eqref{eq:epsilon_disc}--\eqref{eq:eqmu_disc} } \ \forall \, \alpha\in\{1,\dots, N_\mathrm{L}\} \bigg\},
    \label{eq:minprob}
    \end{equation}
where ${i}^\alpha_e\in\{1,\dots, \bar{N}\}$ represents an index pointer $(e,\alpha)\mapsto i^\alpha_e$ that assigns a data point $\bar{\mathbf{z}}_i$ to material point $e$ in loading case $\alpha$. Hereinafter, $\{\Box_{a}^{b}\}$ denotes an array of all elements $\Box$ indexed by $a$ and/or $b$.

Assuming that the mapping index ${i}^\alpha_e\in\{1,\dots, \bar{N}\}$ is given for all loading cases $\alpha\in\{1,\dots,N_\mathrm{L}\}$ and material points $e\in\{1,\dots,M\}$, and imposing the balance equations~\eqref{eq:eqsig_disc}--\eqref{eq:eqmu_disc} via Lagrange multipliers collected in vectors $\discrete{\bm{\lambda}}^{u\alpha}$ and $\discrete{\bm{\lambda}}^{\rchi\alpha}$, the stationarity conditions with respect to $(\sig^\alpha_e,\btau^\alpha_e,\bmu_e^\alpha,\discrete{\bm{\lambda}}^{u\alpha},\discrete{\bm{\lambda}}^{\rchi\alpha})$ read
    \begin{align} 
    &w_e\,\Ctens_e^{-1}\big(\sig^\alpha_e-\bar{\sig}_{i_e^\alpha}\big)-w_e\,\Bea\,\discrete{\bm{\lambda}}^{u\alpha}=\bm{0},\label{eq:sigupd}\\
    &w_e\,\Dtens_e^{-1}\big(\btau^\alpha_e-\bar{\btau}_{i_e^\alpha}\big)-w_e\,\Bga\,\discrete{\bm{\lambda}}^{u\alpha}+w_e\,\Nca\,\discrete{\bm{\lambda}}^{\rchi\alpha}=\bm{0},\label{eq:tauupd}\\
    &w_e\,\Atens_e^{-1}\big(\bmu^\alpha_e-\bar{\bmu}_{i_e^\alpha}\big)-w_e\,\Bza\,\discrete{\bm{\lambda}}^{\rchi\alpha}=\bm{0},\label{eq:muupd}\\       
    &\sum_{e=1}^M w_e\,\Big[\Beta\,\sig_e^\alpha + \Bgta\,\btau_e^\alpha\Big]-\discrete{\mathbf{f}}^\alpha=\bm{0},\label{eq:eqsiga}\\   
    &\sum_{e=1}^M w_e\,\Big[\Bzta\,\bmu_e^\alpha - \Ncta\,\btau_e^\alpha\Big]-\discrete{\mathbf{m}}^\alpha=\bm{0}.\label{eq:eqmua}
    \end{align}           
    Combining equations~\eqref{eq:sigupd}--\eqref{eq:muupd} and~\eqref{eq:eqsiga}--\eqref{eq:eqmua} gives a linear system in $(\discrete{\bm{\lambda}}^{u\alpha},\discrete{\bm{\lambda}}^{\rchi\alpha})$ $\forall\,\alpha\in\{1,\dots,N_\mathrm{L}\}$:
    \begin{align}
    \begin{split}
    \sum_{e=1}^M w_e\,\Big[\Beta\Ctens_e\Bea + \Bgta\Dtens_e\Bga\Big]\discrete{\bm{\lambda}}^{u\alpha} - \sum_{e=1}^M w_e\,\Big[\Bgta\Dtens_e\Nca\Big]\discrete{\bm{\lambda}}^{\rchi\alpha}&\\=\discrete{\mathbf{f}}^\alpha-\sum_{e=1}^M w_e\,\Big[\Beta\bar{\sig}_{i_e^\alpha} + \Bgta\bar{\btau}_{i_e^\alpha}\Big]&,\label{eq:lprob1}\end{split}\\
    \begin{split}
    \sum_{e=1}^M w_e\,\Big[\Ncta\Dtens_e\Bga\Big]\discrete{\bm{\lambda}}^{u\alpha} -\sum_{e=1}^M w_e\,\Big[\Bzta\Atens_e\Bza + \Ncta\Dtens_e\Nca\Big]\discrete{\bm{\lambda}}^{\rchi\alpha}& \\=\sum_{e=1}^M w_e\,\Big[\Bzta\bar{\bmu}_{i_e^\alpha} - \Ncta\bar{\btau}_{i_e^\alpha}\Big] - \discrete{\mathbf{m}}^\alpha&.
    \label{eq:lprob2}\end{split}  
    \end{align}
    The coupled equations~\eqref{eq:lprob1}--\eqref{eq:lprob2} are standard forms analogous to linear micromorphic elasticity with particular right-hand sides, and can thus be easily implemented in conventional finite element programs. The mechanical stresses are then updated locally from~\eqref{eq:sigupd}--\eqref{eq:muupd}.
    
    On the other hand, stationarity with respect to the material data $\{\bar{\mathbf{z}}_i\}$ gives,  $\forall\,i\in\{1,\dots, \bar{N}\}$,
    \begin{align}
        &\sum_{\alpha=1}^{N_\mathrm{L}}\;\sum_{e\in S^\alpha_i} w_e\,\Ctens_e\big(\eps^\alpha_e-\bar{\eps}_{i_e^\alpha}\big)=\bm{0},\qquad
        \sum_{\alpha=1}^{N_\mathrm{L}}\;\sum_{e\in S^\alpha_i} w_e\,\Ctens_e^{-1}\big(\sig^\alpha_e-\bar{\sig}_{i_e^\alpha}\big)=\bm{0},\label{eq:sigstar}\\
        &\sum_{\alpha=1}^{N_\mathrm{L}}\;\sum_{e\in S^\alpha_i} w_e\,\Dtens_e\big(\gam^\alpha_e-\bar{\gam}_{i_e^\alpha}\big)=\bm{0},\qquad\sum_{\alpha=1}^{N_\mathrm{L}}\;\sum_{e\in S^\alpha_i} w_e\,\Dtens_e^{-1}\big(\btau^\alpha_e-\bar{\btau}_{i_e^\alpha}\big)=\bm{0}
        ,\\
        &\sum_{\alpha=1}^{N_\mathrm{L}}\;\sum_{e\in S^\alpha_i} w_e\,\Atens_e\big(\zet^\alpha_e-\bar{\zet}_{i_e^\alpha}\big)=\bm{0},\qquad
        \sum_{\alpha=1}^{N_\mathrm{L}}\;\sum_{e\in S^\alpha_i} w_e\,\Atens_e^{-1}\big(\bmu^\alpha_e-\bar{\bmu}_{i_e^\alpha}\big)=\bm{0},\label{eq:mustar}     
    \end{align}    
    where
    \begin{equation}
        S^\alpha_i = \big\{e\in\{1,\dots,M\} \ \colon \ i^\alpha_e = i \big\}
    \end{equation}
    is a set of material point indices associated with the data point $i$ at snapshot $\alpha$. We thus find from~\eqref{eq:sigstar}--\eqref{eq:mustar} that the $i$th identified material data point $\bar{\mathbf{z}}_i$ is a weighted average of a subset of the mechanically admissible states $\{\mathbf{z}_e^\alpha\}$, restricted to material points and snapshots to which the data point $i$ has been assigned. Moreover, combining equations~\eqref{eq:sigupd}--\eqref{eq:muupd} and~\eqref{eq:sigstar}--\eqref{eq:mustar} gives, $\forall\,i\in\{1,\dots, \bar{N}\}$,
    \begin{align}
    &\sum_{\alpha=1}^{N_\mathrm{L}}\;\sum_{e\in S^\alpha_i} w_e\,\Bea\,\discrete{\bm{\lambda}}^{u\alpha}=\bm{0},\\
    &\sum_{\alpha=1}^{N_\mathrm{L}}\;\sum_{e\in S^\alpha_i} w_e\big(\Bga\,\discrete{\bm{\lambda}}^{u\alpha}-\Nca\,\discrete{\bm{\lambda}}^{\rchi\alpha}\big)=\bm{0},\\
    &\sum_{\alpha=1}^{N_\mathrm{L}}\;\sum_{e\in S^\alpha_i} w_e\,\Bza\,\discrete{\bm{\lambda}}^{\rchi\alpha}=\bm{0}.
    \end{align}
    Accordingly, the Lagrange multipliers are such that the corresponding generalized pseudo-strains have zero $w_e$-weighted average. Finally,  stationarity with respect to $\{i_e^\alpha\}$ corresponds to finding the set of indices assigned to each material point $e$ and load case $\alpha$ such that the distance between the material and mechanical states is minimized. Here, as suggested in~\citet{stainier2019model} for Cauchy continua, minimization with respect to $(\{\bar{\mathbf{z}}_i\},\{i_e^\alpha\})$ is performed simultaneously through $k$-means clustering \tB{using the local distance metric: each mechanical state $\mathbf{z}_e^\alpha$ is assigned by $i_e^\alpha$ to the nearest cluster $i$ in material space, and the identified material data points $\{\bar{\mathbf{z}}_i\}$ follow as the corresponding weighted centroids.}
    
    We present the computational procedure in Algorithm~\ref{alg:fixedpoint}. \tB{Algorithmic convergence of DDI and the influence of the chosen number of data points to be identified, $\bar N$, have been studied recently for the case of Cauchy continua~\cite{hachem2025mathematical,Leygue2025}. Specifically, \citet{hachem2025mathematical} showed for Cauchy continua that a trade-off exists between values of $\bar N$ that are too low, leading to reduced accuracy, and values that are too high, which exacerbate ill-posedness effects in the identification problem. We expect analogous considerations to apply in the present generalized continuum setting.}
    \begin{algorithm}[h!]
    \setstretch{1.35}
    \begin{algorithmic}[1]
    \State Initialize iterations with $j\coloneqq 0$.
    \State Perform $k$-means clustering on $\{(\eps^\alpha_e,\gam^\alpha_e,\zet^\alpha_e)\}$ to find initial data estimates $\{({\bar{\eps}_{i}},{\bar{\gam}_{i}},{\bar{\zet}_{i}})\}^{(0)}$ and $\{i_e^\alpha\}^{(0)}$.
    \State Set an initial guess $\{({\bar{\sig}_{i}},{\bar{\btau}_{i}},{\bar{\bmu}_{i}})\}^{(0)}$.
    \Repeat  
    \State Set $j \leftarrow j+1$.
    \State \multiline{With fixed $\{\bar{\mathbf{z}}_i\}^{(j-1)}$ and $\{i_e^\alpha\}^{(j-1)}$, find $({\discrete{\bm{\lambda}}^{u\alpha}},{\discrete{\bm{\lambda}}^{\rchi\alpha}})^{(j)}$ from~\eqref{eq:lprob1}--\eqref{eq:lprob2} $\forall\,\alpha\in\{1,\dots,N_\mathrm{L}\}$ and update $({\sig_e^\alpha},{\btau_e^\alpha},{\bmu_e^\alpha})^{(j)}$ from~\eqref{eq:sigupd}--\eqref{eq:muupd} $\forall\,\alpha\in\{1,\dots,N_\mathrm{L}\}$ and $e\in\{1,\dots,M\}$.}
    \State Perform $k$-means clustering on $\{\mathbf{z}^\alpha_e\}^{(j)}$ to update the material data $\{\bar{\mathbf{z}}_i\}^{(j)}$ and $\{i_e^\alpha\}^{(j)}$.
    \Until $\{i_e^\alpha\}^{(j)}=\{i_e^\alpha\}^{(j-1)}$.
    \State Set $\{({\sig_e^\alpha},{\btau_e^\alpha},{\bmu_e^\alpha})\}\coloneqq\{({\sig_e^\alpha},{\btau_e^\alpha},{\bmu_e^\alpha})\}^{(j)}$ and $\{\bar{\mathbf{z}}_i\}\coloneqq\{\bar{\mathbf{z}}_i\}^{(j)}$ (identified material data).
    \end{algorithmic}
    \caption{Alternate minimization scheme for problem~\eqref{eq:minprob}.}
    \label{alg:fixedpoint}
    \end{algorithm} 

\section{Numerical simulations}
\label{sec:num}

\tB{This section presents numerical simulations demonstrating the capability of the proposed framework to identify generalized material data from heterogeneous kinematic fields and applied forces. The examples show that, within the class of micromorphic materials encompassed by the local phase space~\eqref{eq:ps} with coordinates $(\eps_e,\gam_e,\zet_e,\sig_e,\btau_e,\bmu_e)$, the method can identify the effective continuum type and the constitutive response associated with its generalized stress--strain components. Specifically, our approach enables faithful identification of non-trivial micromorphic stress features, including non-symmetric and higher-order contributions, across different material behaviors.}

\tB{
\subsection{Synthetic validation}
\label{sec:synthetic}
As a proof of concept, we first conduct synthetic studies using data generated from numerical BVPs governed by ground-truth constitutive models. We aim to assess the ability of the DDI solver to identify the micromorphic constitutive response of both elastic and inelastic materials in a model-agnostic manner, provided that the behavior is described locally within the phase space~\eqref{eq:ps}.

\begin{figure}[b!]
    \centering
    \includeinkscape[width=0.35\textwidth]{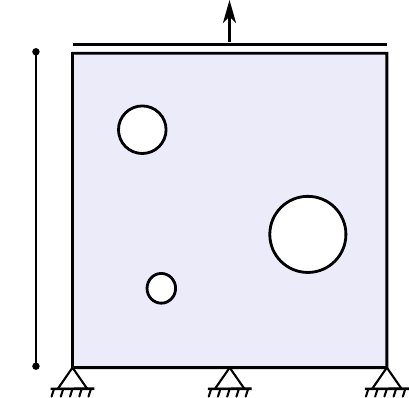}
    \vspace{1em}
    \caption{Boundary value problem for the perforated plate under plane strain conditions. The holes with radii $R_1=0.045H$, $R_2=0.12H$, and $R_3=0.075H$ are located at coordinates $(0.28, 0.25)H$, $(0.72, 0.42)H$, and $(0.22, 0.75)H$, respectively, measured from the bottom-left corner. We assume $H=10$~mm in the simulations.}
    \label{fig:bvp3}
\end{figure}

\subsubsection{Problem setup}
A specimen with geometrical defects is employed, consisting of a plate with three holes under plane strain conditions (Figure~\ref{fig:bvp3}). The sample is meshed using $N=1778$ nodes and 1688 bilinear quadrilateral elements ($M=6752$ material points). For the two examples that follow, we select DDI metric tensors of the form~\eqref{eq:mettens}, with numerical (hyper)parameters chosen regardless of the ground-truth material behavior: $\lambda=86.42$~GPa, $\mu=37.04$~GPa, $c=5$, and $\ell=2/\sqrt{2}$~mm. The number of data points to be identified, $\bar N$, is set to achieve a ratio of mechanical states to material states of~$(MN_\mathrm{L})/\bar{N}\approx100$.

We note that reasonable variations of the metric parameters did not yield significant differences in the examples considered. Automated selection of metric parameters may also be pursued (\ref{sec:metric}), although this approach is not implemented here for simplicity. Similarly, the number of data points $\bar N \in \{2,\dots,MN_\mathrm{L}-1\}$ is chosen heuristically to reduce the data dimensionality by a factor of 100; a systematic study on this topic, in the spirit of~\citet{hachem2025mathematical}, is left for future~work.}

\begin{figure}[b!]
    \centering    \includegraphics[width=0.8\textwidth]{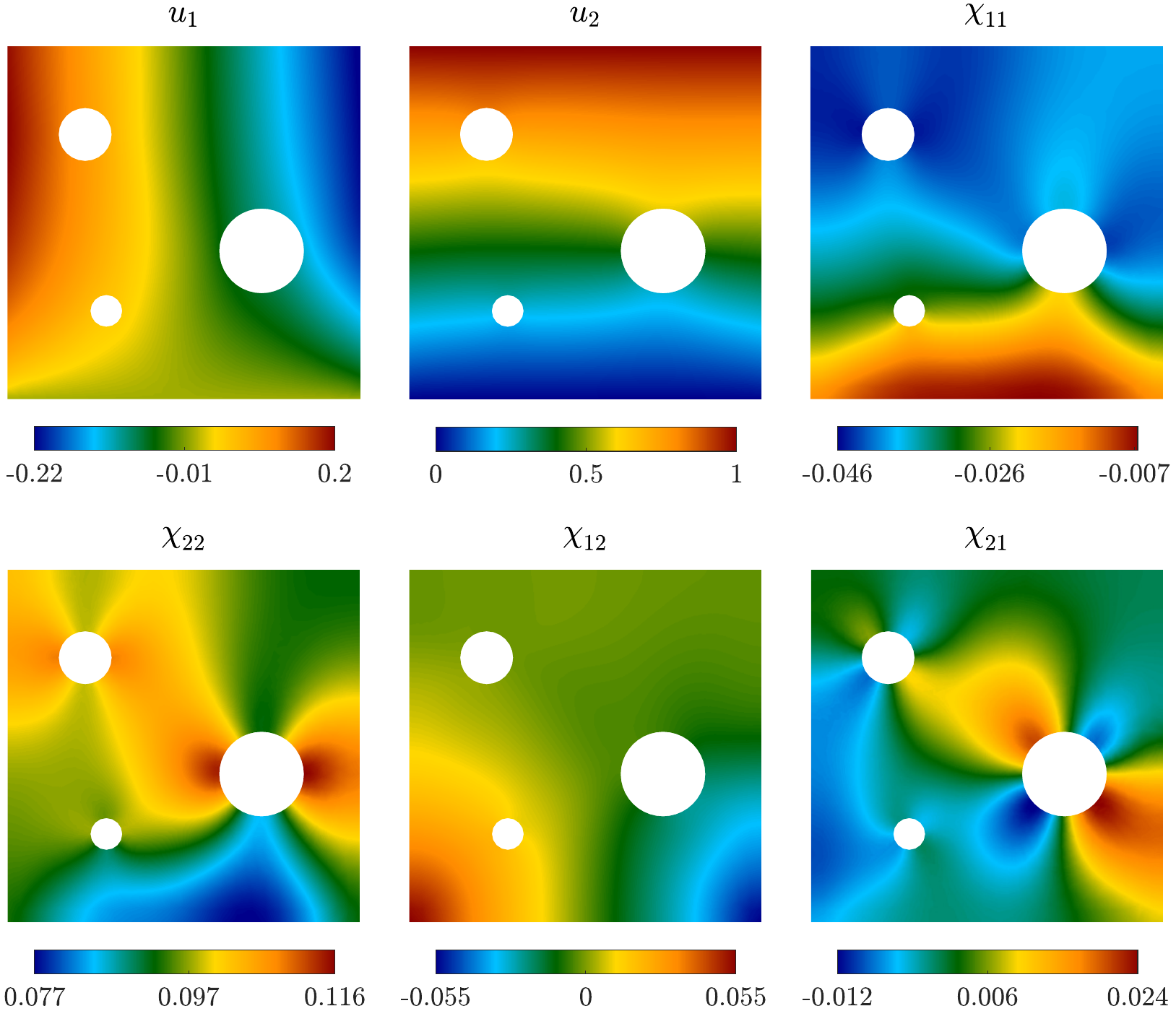}
    \caption{Kinematic data for the perforated plate with a ground-truth elastic micromorphic material, showing the displacements $u_i$ [mm] and non-symmetric micro-deformations $\rchi_{ij}$ at $\bar{u}_2 = 1$~mm.}
    \label{fig:synth_kin_snaps}
\end{figure}

\subsubsection{Elastic micromorphic material}
\label{sec:synthetic_elastic}

\tB{We first consider material behavior obeying a full micromorphic constitutive model, i.e., one that includes both symmetric and skew-symmetric micro-deformations. For the sake of simplicity, and to enable a controlled and systematic quantitative analysis, we assume linear micromorphic elasticity, characterized by a strain energy density
\begin{equation}
    \psi(\eps,\gam,\zet)=\frac{1}{2}\left(\Ctens\,\eps:\eps + \Dtens\,\gam:\gam + \Atens\,\zet\trip\zet\right).
\end{equation}
As such, the constitutive relations naturally follow as
\begin{equation}
    \sig = \frac{\partial\psi}{\partial\eps} = \Ctens\,\eps , \qquad \btau = \frac{\partial\psi}{\partial\gam} = \Dtens\,\gam, \qquad \bmu = \frac{\partial\psi}{\partial\zet} = \Atens\,\zet ,
\end{equation}
providing closure to the governing equations~\eqref{eq:epsilon}--\eqref{eq:eqmu}. We assume ground-truth constitutive tensors
\begin{subequations}\label{eq:eltens}
\begin{align}
     &\mathsf{C}_{ijkl} \;=\; \lambda \, \delta_{ij}\delta_{kl} \;+\; \mu\, (\delta_{ik}\delta_{jl} + \delta_{il}\delta_{jk}),   
     \\ 
     &\mathsf{D}_{ijkl} \;=\; c_1\,\mathsf{C}_{ijkl} \;+\; c_2\,(\delta_{ik}\delta_{jl} - \delta_{il}\delta_{jk}),
     \\
     &\mathsf{A}_{ijklmn} \;=\; \ell_1^2\,\mathsf{C}_{ijlm}\,\delta_{kn} \;+\; \mu\,\ell_2^2 \, (\delta_{il}\delta_{jm} - \delta_{im}\delta_{jl})\,\delta_{kn},
\end{align}
\end{subequations}
with Young's modulus $E=217.5$~GPa and Poisson's ratio $\nu=0.3$ (i.e., $\lambda=125.48$~GPa, $\mu=83.65$~GPa). We further set $c_1=4.26$, $c_2=356.63$~GPa, and $\ell_1=\ell_2=2/\sqrt{2}$~mm.}

\tB{To generate synthetic kinematic and boundary force data from the model above,} finite element simulations are conducted by imposing displacements $\bar{u}_2$ in $N_\mathrm{L}=30$ increments from 0 to 1~mm. The micro-deformations are left free on the boundary. At each load step, snapshots of the kinematic fields are recorded as nodal displacements and micro-deformations $(\discrete{\mathbf{u}},\discrete{\bm{\rchi}})$, together with the corresponding reaction forces on the top side. Figure~\ref{fig:synth_kin_snaps} shows a snapshot of the kinematics at $\bar{u}_2=1$~mm. \tB{Note, in particular, the non-symmetric micro-deformation components $\rchi_{ij}$ arising from the polar coupling modulus $c_2$.} The generalized mechanical strains $(\eps_e,\gam_e,\zet_e)$ then follow locally from~\eqref{eq:epsilon_disc}--\eqref{eq:zeta_disc}.

\begin{figure}[b!]
    \centering
    \begin{subfigure}[t]{0.45\textwidth}
        \includegraphics[width=\textwidth]{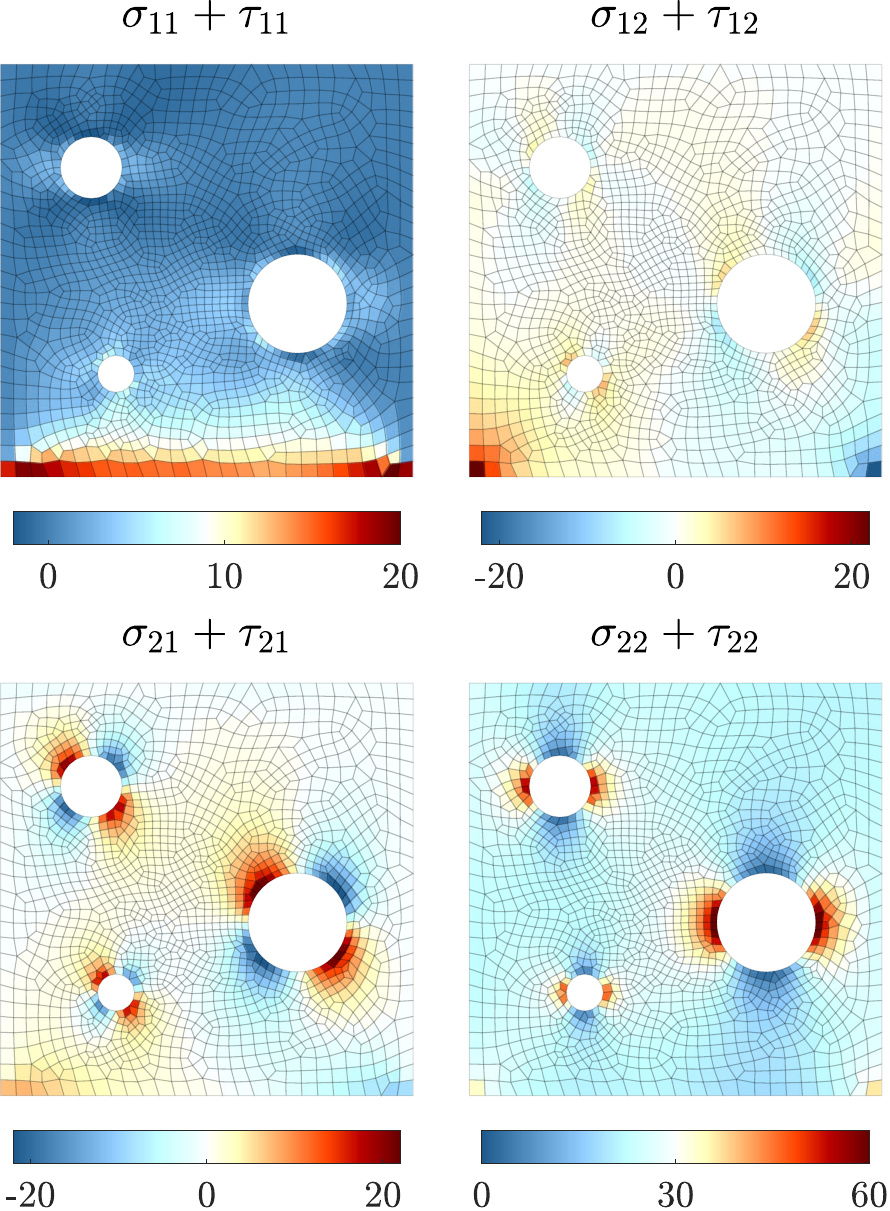}
        \subcaption{Reference}
    \end{subfigure}
    \hspace{0.5em}
    \vline
    \hspace{1em}
    \begin{subfigure}[t]{0.45\textwidth}
        \includegraphics[width=\textwidth]{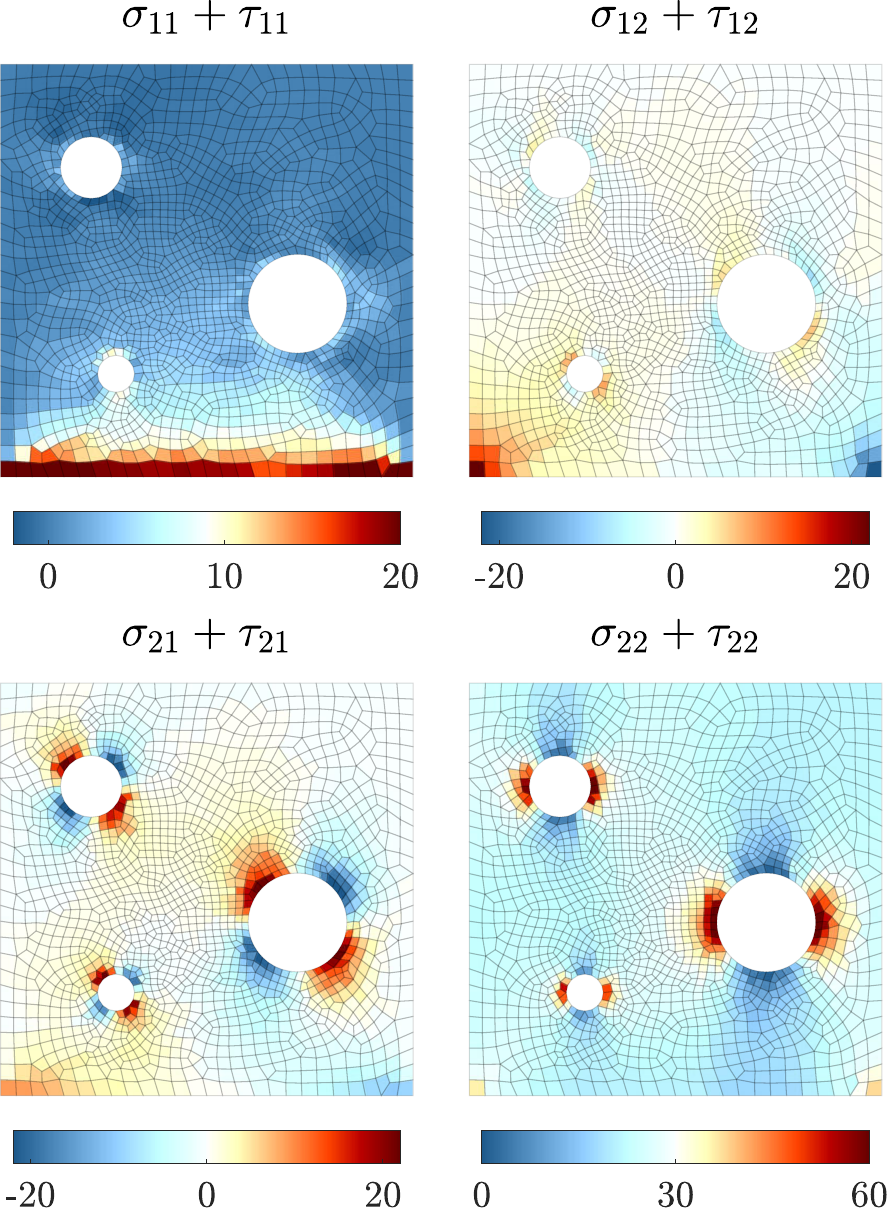}
        \subcaption{Identified}
    \end{subfigure}
    \caption{Comparison of generalized second-order stresses $\sigma_{ij}+\tau_{ij}$ [GPa], showing (a) the ground-truth reference simulation and (b) the identified stress fields.}
    \label{fig:synth_s_snaps}
\end{figure}

\begin{figure}[h!]
    \centering
    \begin{subfigure}[t]{0.9\textwidth}
        \includegraphics[width=\textwidth]{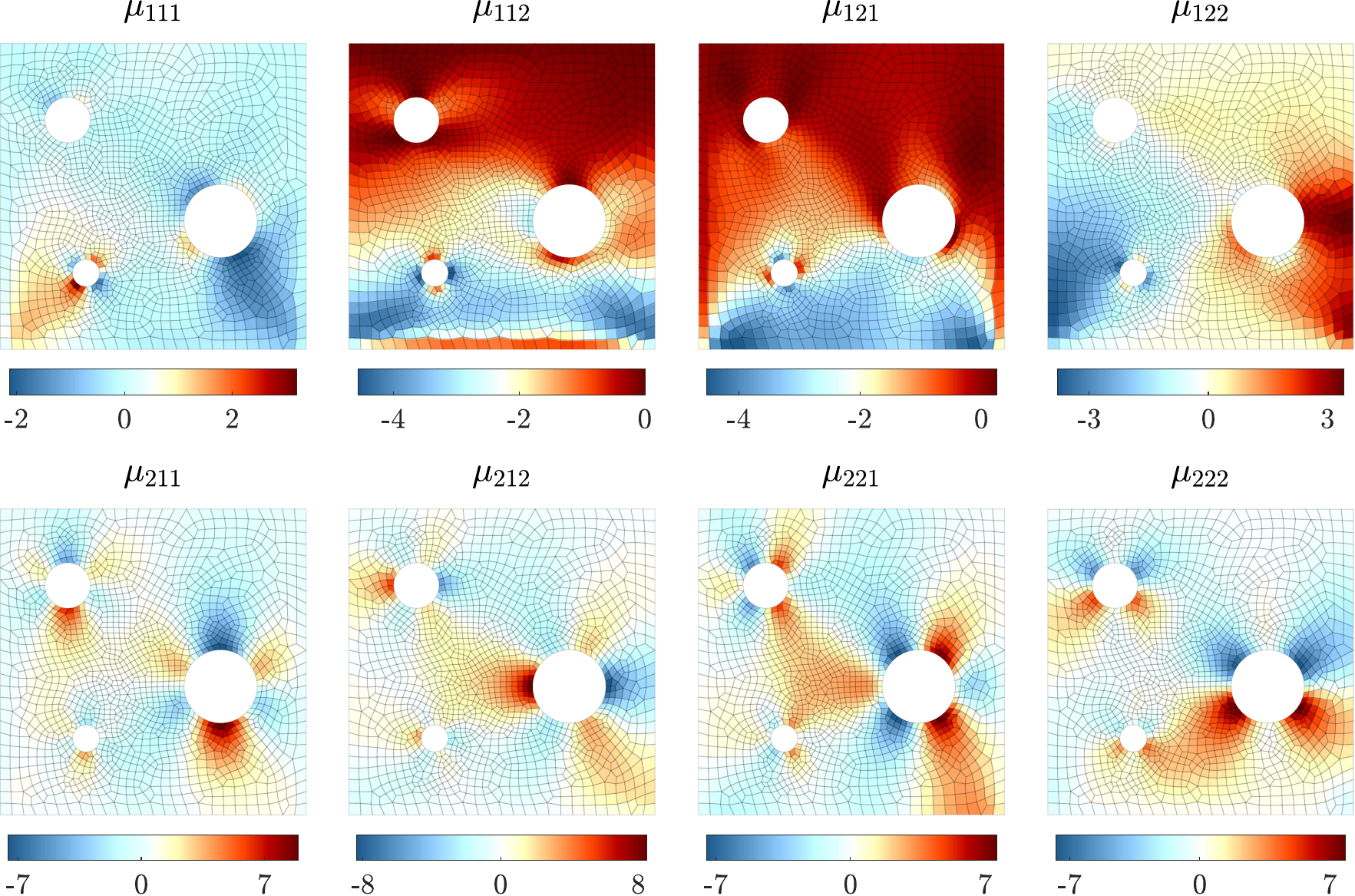}
        \subcaption{Reference}
    \end{subfigure}
    \vskip 0.5em
   \noindent\rule{0.9\textwidth}{0.4pt}
   \vskip 1em
    \begin{subfigure}[t]{0.9\textwidth}
        \includegraphics[width=\textwidth]{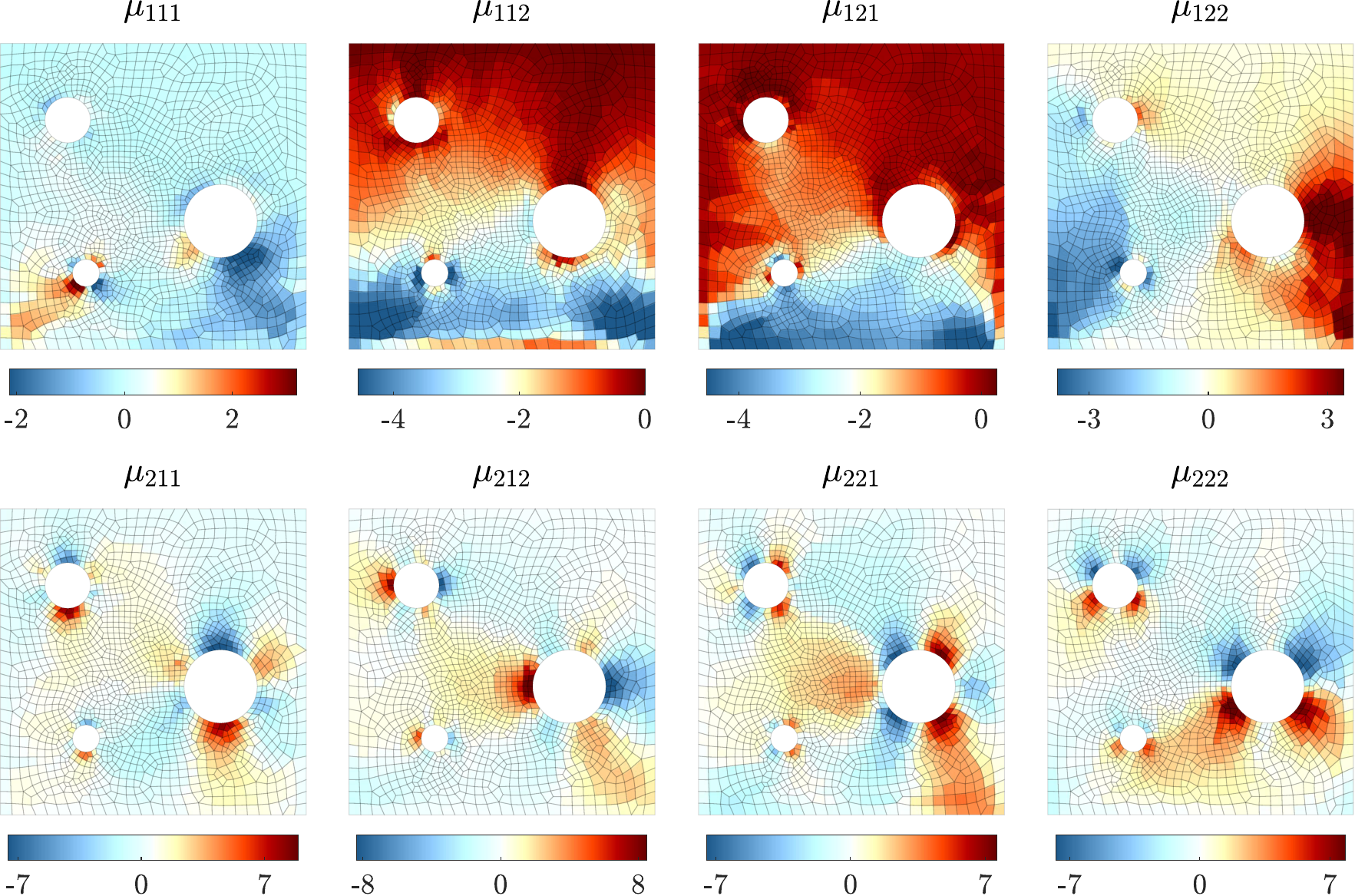}
        \subcaption{Identified}
    \end{subfigure}
    \caption{Comparison of generalized third-order stresses $\mu_{ijk}$ [GPa$\,\cdot\,$mm], showing (a) the ground-truth reference simulation and (b) the identified double stress fields.}
    \label{fig:synth_mu_snaps}
\end{figure}

\tB{This simulated kinematic data and the corresponding boundary forces are passed to the DDI solver to obtain the corresponding stress fields and build a representative material dataset of $\bar{N}=2026$ points.} Figure~\ref{fig:synth_s_snaps} shows the generalized second-order stress fields $\sigma_{ij} + \tau_{ij}$ at $\bar{u}_2=1$~mm, comparing the identified mechanical states to the reference constitutive model. While a very slight noise level is observed, the contour plots show remarkable agreement between the identified and reference fields. In particular, the non-symmetric nature of $\tau_{ij}$ is properly identified. Moreover, although an overestimation of $\sigma_{11}+\tau_{11}$ is present along the bottom edge, the stress concentrations show good agreement in both location and magnitude. Figure~\ref{fig:synth_mu_snaps} presents the generalized third-order stress fields $\mu_{ijk}$, again contrasting the identified mechanical states with the reference constitutive response. The level of noise in the identified response is higher than for the second-order components. Nevertheless, excellent agreement is still observed, both in terms of spatial distribution and magnitude of the different third-order components.

\begin{figure}[t!]
    \centering
        \begin{overpic}[width=0.95\textwidth]{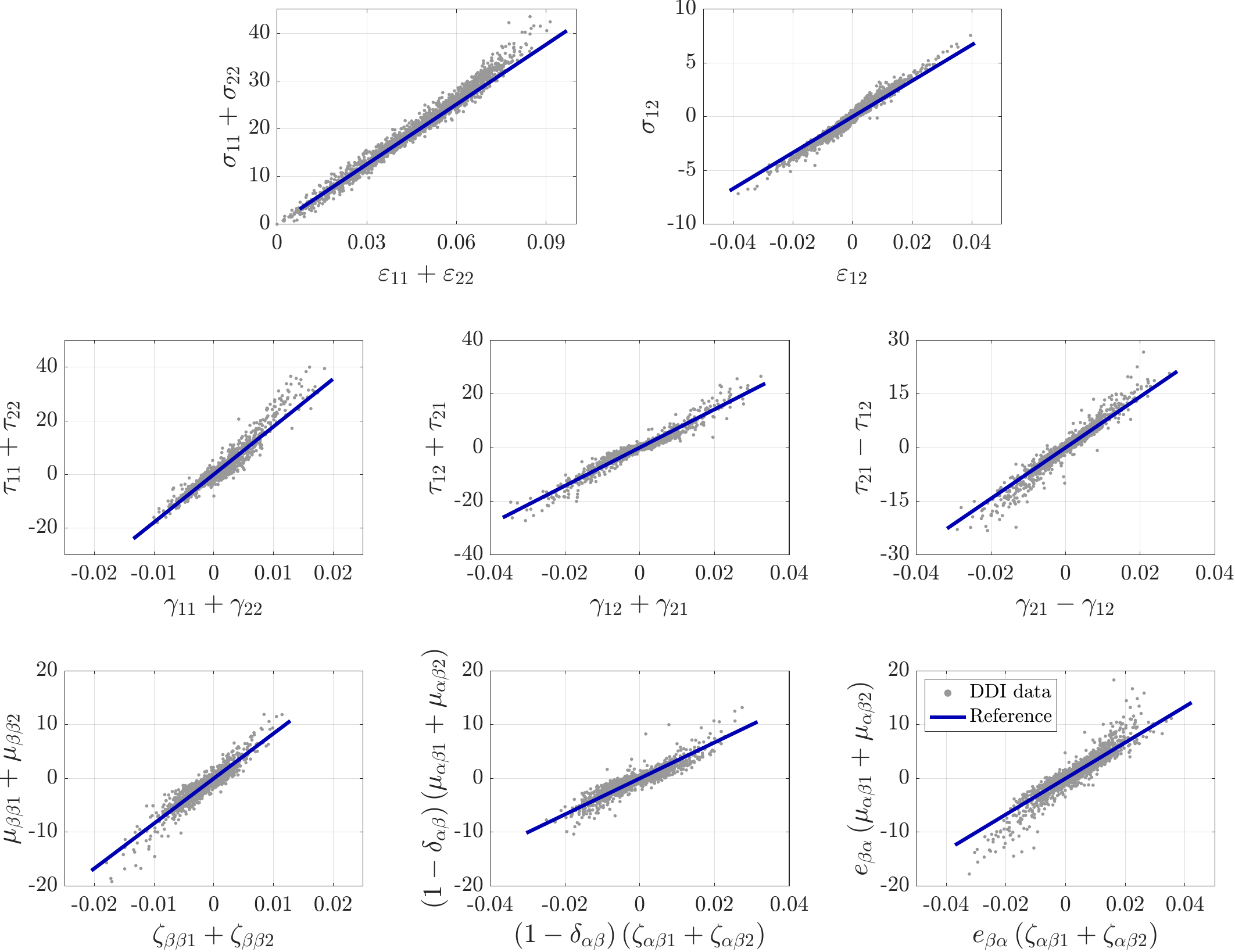}
        \put(0,-0.5){\small(c1)}
        \put(33,-0.5){\small(c2)}
        \put(66,-0.5){\small(c3)}
        \put(0,27){\small(b1)}
        \put(33,27){\small(b2)}
        \put(66,27){\small(b3)}
        \put(18,54){\small(a1)}
        \put(52,54){\small(a2)}
        \end{overpic}
        \caption{Material data identified from all snapshots of the perforated plate (gray, DDI) and expected behavior from the linear micromorphic ground-truth constitutive model (blue). Note that summation over repeated indices is assumed, where $\alpha,\beta\in\{1,2\}$. Stress quantities are shown in GPa. See Table~\ref{tab:synthDDI_all_scatter} for representative statistics.}
    \label{fig:synthDDI_all_scatter}
\end{figure}

\begin{table}[b!]
\caption{Linear relations of the reference elastic micromorphic model and summary statistics of the identified data: $R^2$ values from least-squares regression and relative errors with respect to the ground-truth constitutive parameters.}
\centering
\renewcommand{\arraystretch}{1.1}
\begin{tabular}{lllcc}
\toprule
Relation & Fig.~\ref{fig:synthDDI_all_scatter} & Ref.~slope & $R^2$ & Slope error [\%] \\
\hline
$\vphantom{\dfrac{\sum}{\sum}}(\sigma_{11}+\sigma_{22})/(\varepsilon_{11}+\varepsilon_{22})$
    & (a1) & $2(\lambda + \mu)$ & 0.98 & 3.53 \\[4pt]

$\sigma_{12}/\varepsilon_{12}$
    & (a2) & $2\mu$ & 0.97 & 19.06 \\[4pt]

$(\tau_{11}+\tau_{22})/(\gamma_{11}+\gamma_{22})$
    & (b1) & $2c_1(\lambda + \mu)$ &  0.91 & 6.34 \\[4pt]

$(\tau_{12}+\tau_{21})/(\gamma_{12}+\gamma_{21})$
    & (b2) & $2c_1\mu$ & 0.93 & 3.48 \\[4pt]

$(\tau_{21}-\tau_{12})/(\gamma_{21}-\gamma_{12})$
    & (b3) & $2c_2$ &  0.93 & 13.67 \\[4pt]

$(\mu_{111}+\mu_{112}+\mu_{221}+\mu_{222})/(\zeta_{111}+\zeta_{112}+\zeta_{221}+\zeta_{222})$
    & (c1) & $2\ell_1^2(\lambda + \mu)$ & 0.85 & 0.98 \\[4pt]

$(\mu_{121}+\mu_{122}+\mu_{211}+\mu_{212})/(\zeta_{121}+\zeta_{122}+\zeta_{211}+\zeta_{212})$
    &  (c2) & $2\mu\,\ell_1^2$ & 0.85 & 8.42 \\[4pt]

$(\mu_{211}+\mu_{212}-\mu_{121}-\mu_{122})/(\zeta_{211}+\zeta_{212}-\zeta_{121}-\zeta_{122})$
    & (c3) & $2\mu\,\ell_2^2$ & 0.86 & 10.72 \\[4pt]
\bottomrule
\end{tabular}
\label{tab:synthDDI_all_scatter}
\end{table}

Recall that in the reference constitutive model, the volumetric and deviatoric parts of the classical stress--strain relation $\sig = \Ctens\,\eps$; the volumetric, symmetric–deviatoric, and skew-symmetric parts of the relative stress--strain relation $\btau = \Dtens\,\gam$; and the volumetric, symmetric–deviatoric, and skew-symmetric gradient parts of the higher-order relation $\bmu = \Atens\,\zet$ (with the decomposition taken over the first two indices) are all linear. Figure~\ref{fig:synthDDI_all_scatter} shows representative point clouds on the associated phase-space planes, considering all generalized stress--strain components. Specifically, we plot the identified data points $\{\bar{\mathbf{z}}_i\}$ against the expected linear behavior from the ground-truth model. The identified material data aligns closely with the reference behavior. Nevertheless, some spread is observed, attributed to the heterogeneity of the input kinematic fields (cf.~\citet{stainier2019model}). 

Table~\ref{tab:synthDDI_all_scatter} reports the correlation values for each linear pair, as well as the relative error in the corresponding expected slopes. The high $R^2$ values indicate that the linear behavior is well captured. The predicted slopes, obtained via least-squares fitting of the data points, closely approximate the expected values from the linear relations, with the highest relative error observed for the $(\varepsilon_{12},\sigma_{12})$ components. Specifically, as visible in Figure~\ref{fig:synthDDI_all_scatter}a2, the corresponding data points align more steeply at small shear stress/strain values, leading to an overestimation of the slope in a least-squares fit. A similar effect, albeit to a lesser degree, is observed for the skew-symmetric components in figures~\ref{fig:synthDDI_all_scatter}b3 and~\ref{fig:synthDDI_all_scatter}c3, while all the other components show a much closer slope agreement. These observations notwithstanding, we emphasize that the objective of the proposed framework is to remain model-free; accordingly, we aim to identify relevant data points rather than material parameters. Indeed, the method remains applicable in practical scenarios where the ground-truth constitutive model is unknown and possibly nonlinear.

\FloatBarrier

\tB{
\subsubsection{Elastoplastic microstrain material}
\label{sec:synthetic_plastic}

\begin{figure}[b!]
    \centering
    \includegraphics[width=0.4\textwidth]{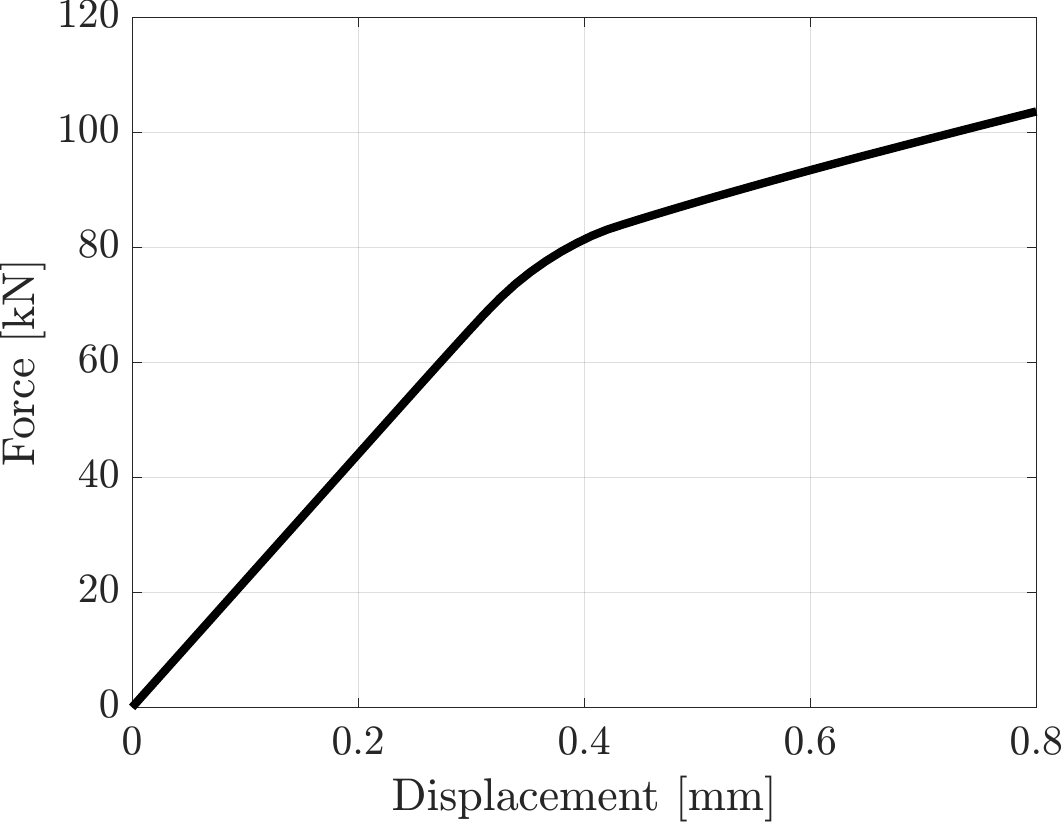}
    \caption{Global force--displacement curve from the perforated plate with a ground-truth elastoplastic microstrain material.}
    \label{fig:res_plastic_fd}
\end{figure}

We now consider nonlinear material behavior obeying a microstrain elastoplastic model, for which the relative strain is symmetric, $\gamma_{ij}=\gamma_{ji}$, and the micro-deformation gradient inherits the same symmetry in its first two indices, $\zeta_{ijk}=\zeta_{jik}$. To introduce inelasticity, we assume that the classical strain admits a plastic contribution, such that $\eps=\eps^\mathrm{e}+\eps^\mathrm{p}$, whereas the relative and higher-order strains remain elastic. Models of this type have been considered, e.g., by~\citet{dillard2006} in the context of metallic foams. The free energy of the model is given by
\begin{equation}
    \psi(\eps,\gam,\zet,\eps^\mathrm{p},p)
    =\frac{1}{2}\Big(
    \Ctens\,(\eps-\eps^\mathrm{p}):(\eps-\eps^\mathrm{p})
    +\Dtens\,\gam:\gam
    +\Atens\,\zet\trip\zet
    + Hp^2 \Big),
\end{equation}
where $H$ is the hardening modulus and $p$ is the equivalent plastic strain. Hence, $\sig=\Ctens\,(\eps-\eps^\mathrm{p})$, $\btau=\Dtens\,\gam$, and $\bmu=\Atens\,\zet$. Plastic flow is governed by the von-Mises yield function
\begin{equation}
    f(\sig;p)=\sqrt{\frac{3}{2}}\Vert\mathrm{dev}\,\sig\Vert-\big(\sigma^\mathrm{p}+Hp\big)\leq0,
\end{equation}
where $\sigma^\mathrm{p}$ is the initial yield strength. The evolution of the plastic strain tensor $\eps^\mathrm{p}$ is assumed associative and rate-independent. We further assume constitutive tensors of the form~\eqref{eq:eltens}, setting $c_2=0$ and $\ell_2=0$ to deactivate the micropolar parts and otherwise employing the same elasticity properties as in the previous example. The plasticity parameters are set to $\sigma^\mathrm{p}=7500$~MPa and $H=30$~GPa.

To generate the synthetic kinematic and boundary force data, finite element simulations are conducted imposing displacements $\bar{u}_2$ in $N_\mathrm{L}=60$ increments from 0 to 0.8~mm. Figure~\ref{fig:res_plastic_fd} shows the resulting nonlinear global force--displacement curve, with the force taken as the reaction resultant on the top side.

\begin{figure}[b!]
    \centering
    \includegraphics[width=0.95\textwidth]{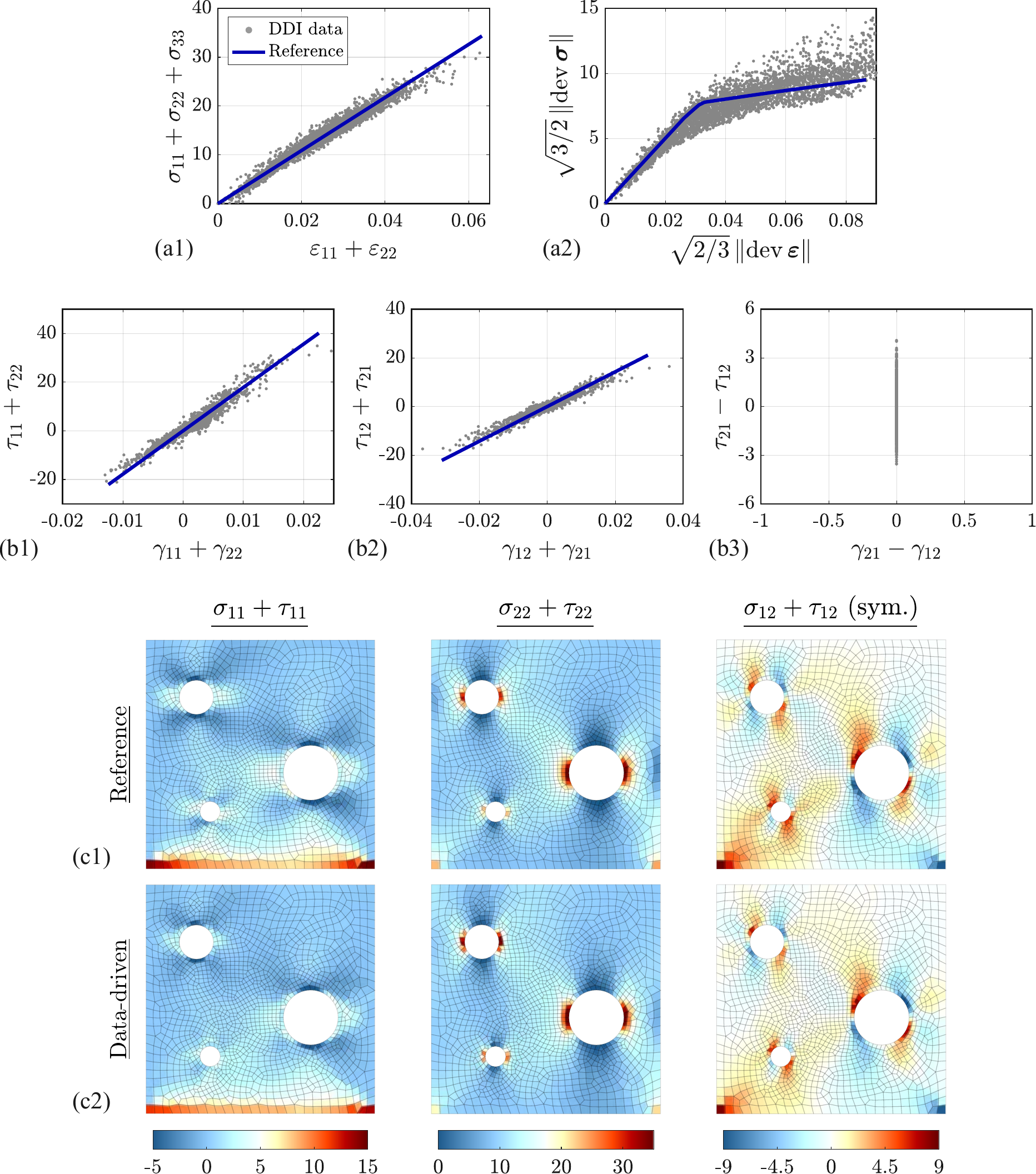}
    \caption{\tB{(a--b) Second-order material data identified from all snapshots of the perforated plate (gray, DDI) and expected behavior from the elastoplastic microstrain ground-truth constitutive model (blue). (c) 
    Comparison of the generalized second-order stress fields $\sigma_{ij}+\tau_{ij}$ between the ground-truth reference simulation and the identified mechanical states. Stress quantities are shown in GPa.}}
    \label{fig:res_plastic_s}
\end{figure}

The simulated kinematic data and the corresponding boundary forces are passed to the DDI solver to obtain the generalized stress fields and build a representative material dataset, using the same algorithmic choices and hyperparameters as in the previous example. Figure~\ref{fig:res_plastic_s}a--b shows selected phase-space projections of the identified second-order material data compared with the reference response. Some noise notwithstanding, the identified data recovers the nonlinear classical deviatoric stress--strain response, including the transition from elastic to plastic behavior and the hardening slope  (Figure~\ref{fig:res_plastic_s}a). Moreover, the relative stress--strain data remains consistent with the elastic microstrain response in the active symmetric components (Figure~\ref{fig:res_plastic_s}b1--2), while the inactive skew-symmetric/micropolar component exhibits a vertical alignment consistent with microstrain kinematics~{(Figure~\ref{fig:res_plastic_s}b3)}. The corresponding second-order stress fields also agree closely with the reference solution (Figure~\ref{fig:res_plastic_s}c). Figure~\ref{fig:res_plastic_mu} further shows the higher-order response. The identified higher-order material data follows the expected microstrain behavior (Figure~\ref{fig:res_plastic_mu}a), and the corresponding third-order stress fields $\mu_{ijk}$ are recovered with reasonable agreement in spatial distribution and magnitude (Figure~\ref{fig:res_plastic_mu}b). We conclude from these results that, while operating in the full micromorphic phase space with a standard metric, the proposed DDI framework can identify nonlinear constitutive behavior in active components while inactive components remain irrelevant.

\begin{figure}[t!]
    \centering
    \includegraphics[width=\textwidth]{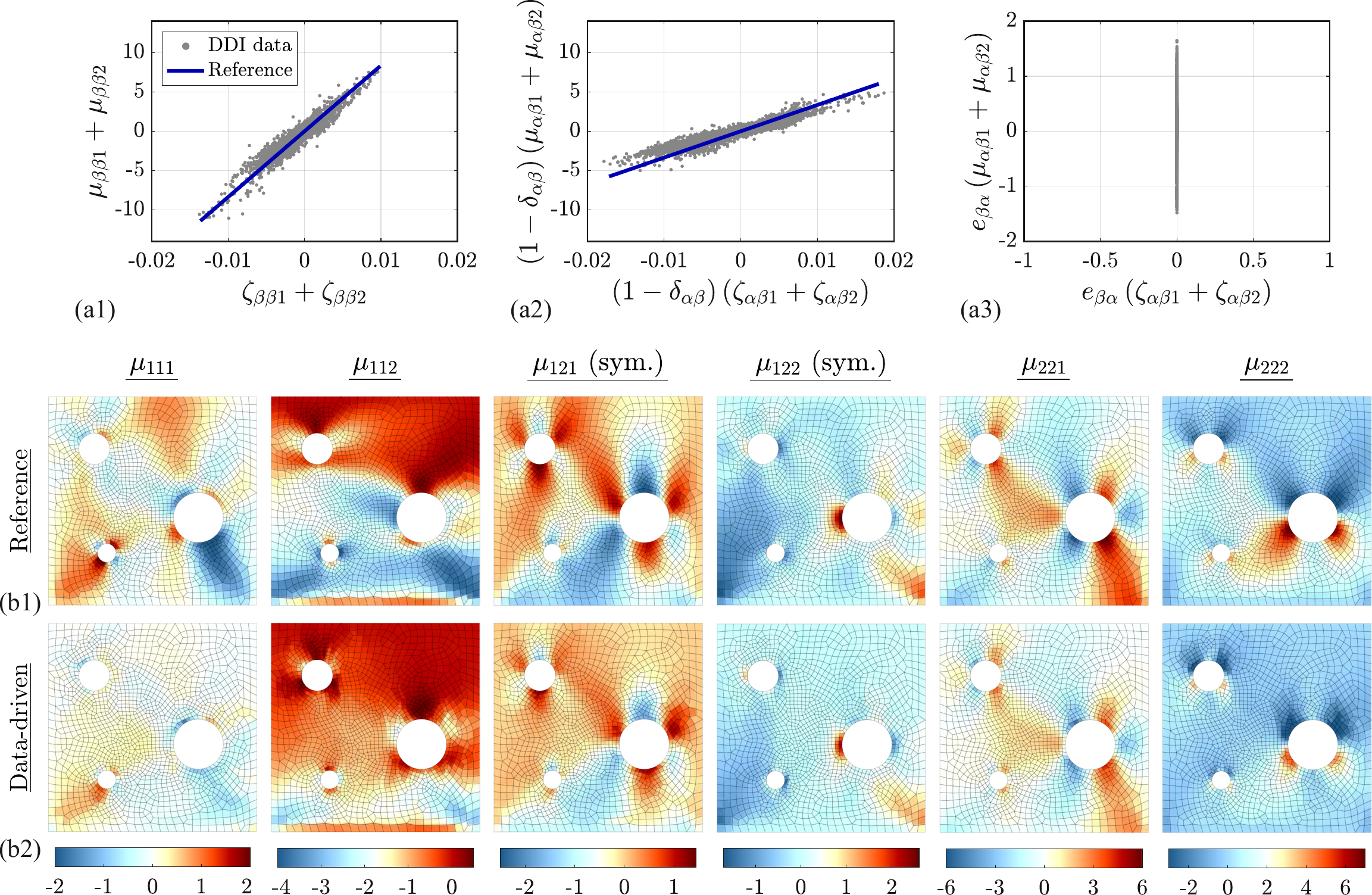}
    \caption{\tB{(a) Higher-order material data identified from all snapshots of the perforated plate (gray, DDI) and expected microstrain behavior from the ground-truth constitutive model (blue). (b) Comparison of the generalized third-order stress fields $\mu_{ijk}$ between the ground-truth reference simulation and the identified mechanical states, with symmetry taken on the first two indices. Double-stress quantities are shown in GPa$\,\cdot\,$mm.}}
    \label{fig:res_plastic_mu}
\end{figure}
}

\FloatBarrier

\subsection{Identification and prediction of mechanical metamaterials}
\label{sec:MMM}
As an example of application, the proposed framework is now used to identify and predict the behavior of mechanical metamaterials in a reduced phase space. Specifically, we consider BVPs where the material is composed of an elastic honeycomb lattice (Figure~\ref{fig:bvp2}). Suppose the base lattice material obeys Euler--Bernoulli beam theory. Then, zeroth-order variational homogenization~\cite{ariza2024homogenization} yields an explicit effective strain energy density that \tB{serves here as a benchmark} for the continuum behavior of the metastructure:
\begin{equation}
    \psi(\eps,\gam) 
     = \frac{1}{2}\left( \Ctens\,\eps:\eps \;+\; \kappa\,\mathrm{skw}\gam:\mathrm{skw}\gam \right)
            = \frac{1}{2} \Ctens\,\eps:\eps \;+\; \kappa\big((\mathrm{curl}\,\bm{u})_3/2-\theta_3\big)^2,
    \label{eq:psi0}
\end{equation}
where, in 2D, $(\mathrm{skw}\,\gam)_{21} = (\mathrm{curl}\,\bm{u})_3/2 - \theta_3$, with $(\mathrm{curl}\,\bm{u})_3 = u_{2,1} - u_{1,2}$, while $\theta_3 = (\rchi_{21} - \rchi_{12})/2$ is an out-of-plane micro-rotation field. The effective moduli are given in Voigt form~by~\cite{ariza2024homogenization}
\begin{equation}
    C_{11}
    =
    \frac
    {
        {EA} \, \left(\, {EA} \, L^2 + 36 \, {EI} \, \right)
    }
    {
        2 \sqrt{3} \, \left(\, {EA} \, L^3 + 12 \, {EI} \, L\right)
    } ,
    \quad
    C_{12}
    =
    \frac
    {
        {EA} \, \left({EA} \, L^2 - 12 \, {EI} \, \right)
    }
    {
        2 \sqrt{3} \, 
        \left(\, {EA} \, L^3 + 12 \, {EI} \, L\right)
    } ,
    \quad
    C_{33}
    =
    \frac
    {
        4 \sqrt{3} \, {EA} \, {EI}
    }
    {
        {EA} \, L^3 + 12 \, {EI} \, L
    } ,
    \label{eq:moduli0}
\end{equation}
while the polar modulus reads
\begin{equation}
        \kappa =  \frac{4\sqrt{3} EI}{L^3},
    \label{eq:moduli1}
\end{equation}
where $E$, $A$, $I$, and $L$ denote the Young's modulus, cross-sectional area, moment of inertia, and length of the beams comprising the microstructure. This homogenized model $\Gamma$-converges to the behavior of the discrete metastructure in the continuum limit $\epsilon\to0$, where $\epsilon$ denotes the ratio between the unit-cell size and the characteristic size of the macroscopic domain~\cite{ariza2024homogenization,Ulloa2025}.

\tB{A key feature of this material response is that non-compatible micro-deformations are purely rotational, arising physically from joint rotations in the lattice. The relative strains are thus purely skew-symmetric,
\begin{equation}
\gamma_{\alpha\beta} = (\mathrm{skw}\,\nabla\bm u)_{\alpha\beta}-e_{\beta\alpha}\,\theta_3 \quad \text{(2D)}, \qquad \gamma_{ij} = (\mathrm{skw}\,\nabla\bm u)_{ij} - e_{jik}\,\theta_k \quad \text{(3D)},
\label{eq:microrot}
\end{equation}
where $e_{\alpha\beta}$ and $e_{ijk}$ denote the components of the Levi-Civita permutation tensor in 2D and 3D, respectively.}  Moreover, the energy~\eqref{eq:psi0} corresponds to a particular generalized continuum of order zero\md a micropolar material independent of curvature effects, and thus insensitive to the length-scale $\epsilon$ (see~\citet{Ulloa2024b} for higher-order generalizations of the theory and size effects). 

\tB{This explicit analytical setting supplied by homogenization theory is convenient for validating the proposed identification framework, since the continuum energy~\eqref{eq:psi0} and coefficients~\eqref{eq:moduli0}--\eqref{eq:moduli1} are known in closed form for the honeycomb lattice and can be used to generate benchmark solutions. Here, the generalized phase-space coordinates are given by the quadruple $(\eps, \gam, \sig, \btau)$. For more general lattice geometries, this effective continuum form still applies, but the explicit energy function and associated coefficients may not be available. In such cases, the data-driven identification process remains applicable, as long as kinematic data and the corresponding external forces are available from BVPs defined on the metastructure. 

Next, for the case of honeycomb lattices, we identify material data from a BVP in the reduced $(\eps, \gam, \sig, \btau)$ phase space, use the identified dataset to predict the response of another, unseen BVP, and systematically compare both identification and prediction results to the corresponding reference solutions.}

\subsubsection{Data identification}

\begin{figure}[t!]
    \centering
    \includeinkscape[width=0.9\textwidth]{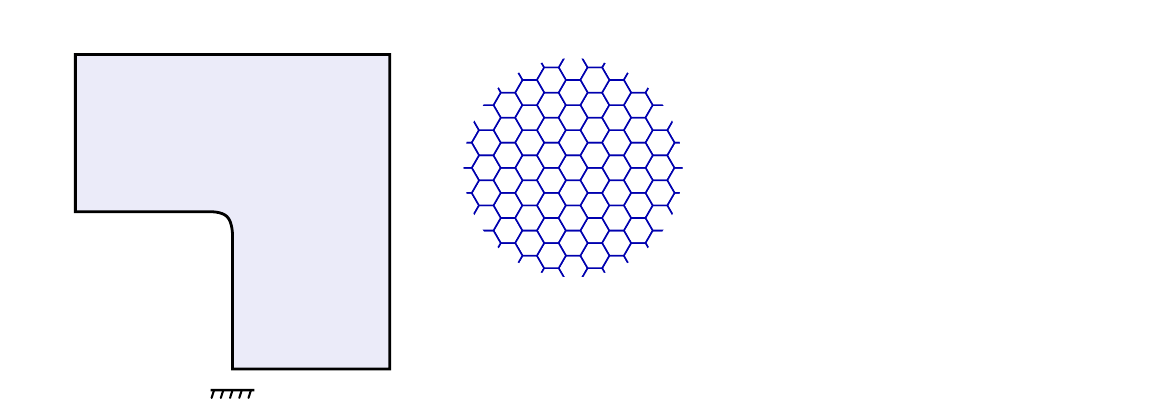}
    \begin{picture}(0,0)
        \put(-395,0){(a)}
        \put(-150,0){(b)}
    \end{picture}    
    \vspace{1em}
    \caption{Schematic representation of boundary value problems with a honeycomb microstructure: (a) L-shaped specimen used for material identification, and (b) double-notched specimen used for data-driven predictions. We assume $H = 30$~mm and $R = 6$~mm. The L-shaped specimen includes a small fillet of radius $0.75$~mm at the inner corner.}
    \label{fig:bvp2}
\end{figure}

\begin{figure}[t!]
    \centering
    \includegraphics[width=0.8\textwidth]{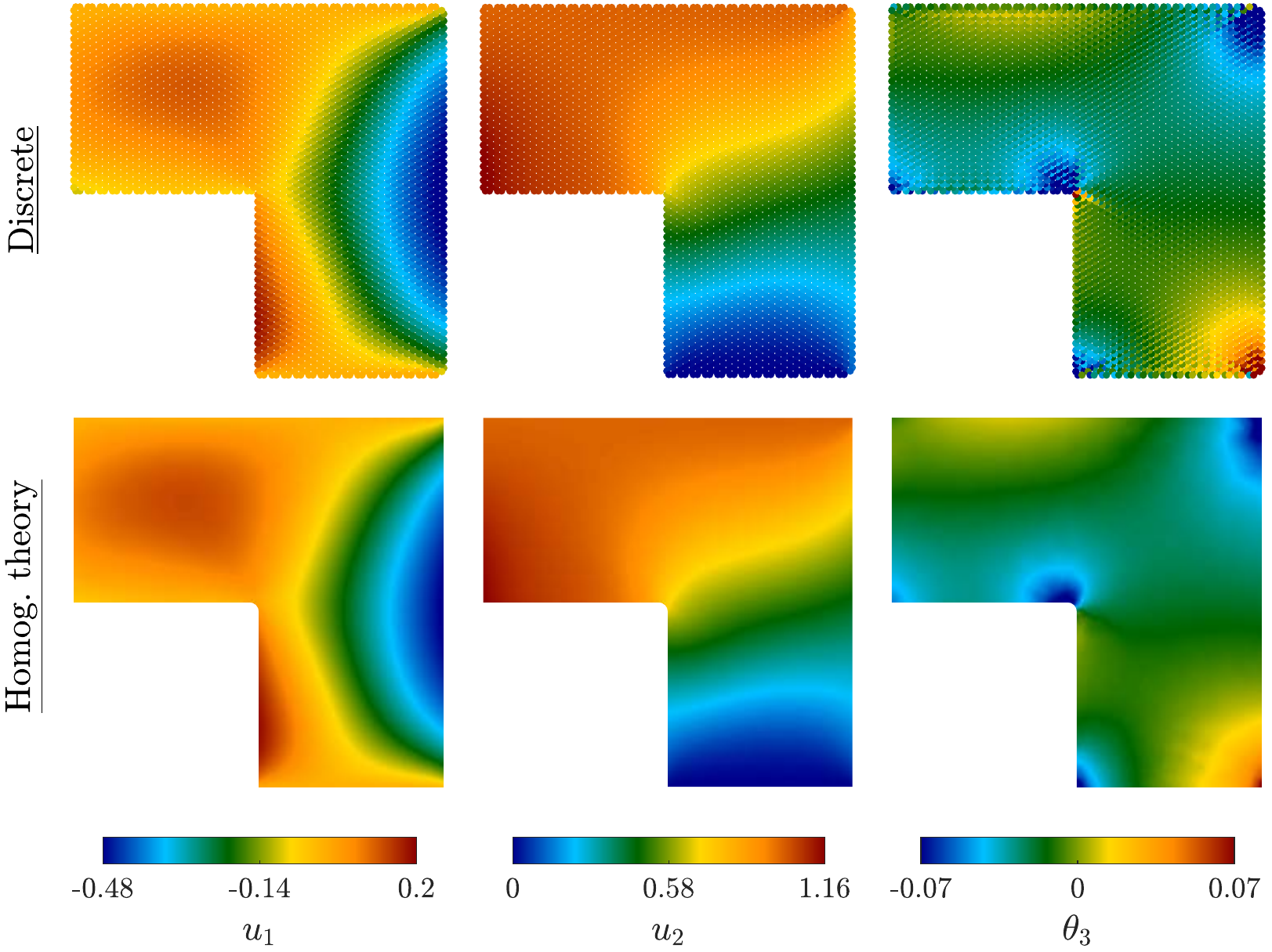}
    \caption{Kinematic data: components $u_i$ [mm] and micro-rotations $\theta_3$ of the honeycomb L-shaped specimen at $\bar{u}_2=1$~mm, showing the direct numerical simulation of the discrete metastructure (top) and the homogenized model (bottom).}
    \label{fig:honey_kin_snaps}
\end{figure}

Consider the L-shaped specimen shown in Figure~\ref{fig:bvp2}a. At the lattice scale, the honeycomb metamaterial is comprised of straight beams with Young's modulus $E=430$~MPa, cross-sectional area $A=0.2$~mm$^2$, second moment of inertia $I=6.67\times10^{-4}$~mm$^4$, and length $L=2$~mm, leading to homogenized model parameters $\lambda=12.167$~MPa, $\mu=0.246$~MPa, and $\kappa=0.248$~MPa  (equations~\eqref{eq:psi0}--\eqref{eq:moduli1}). The bottom side is fixed in both directions, while vertical displacements $\bar{u}_2$ are applied to the top side in $N_\mathrm{L}=30$ increments from 0 to 1~mm. Moreover, a uniformly distributed body couple $\bm{c}=C\bm{e}_3$ is applied in the out-of-plane direction, such that, in equation~\eqref{eq:eqmu},
\begin{equation}
M_{ij} = -e_{ijk} \, c_{k}.
\end{equation}
This body couple could arise, for example, from the action of a horizontal magnetic field $\bm{B}=B\bm{e}_1$ on vertically aligned dipoles $\bm{m}=M\bm{e}_2$ distributed over the joints of the lattice. In the honeycomb geometry, each unit cell contains two node types~\cite{ariza2024homogenization}. Denoting the unit-cell volume by $V$, the resulting body couple per unit volume is
\begin{equation}
\bm{c} = \frac{2}{V}\,\bm{m}\times\bm{B} = -\frac{2}{V} MB\,\bm{e}_3, \qquad V = \frac{3\sqrt{3}}{2}L^2.
\end{equation}
Here, we set $MB = 0.03$~N$\,\cdot\,$mm, corresponding to a distributed couple magnitude~$C=-0.0058$~N/mm$^2$.

Recall that the discrete model, with domain size and applied displacement scaled by $\epsilon$, must converge to the homogenized model response as $\epsilon \to 0$. Figure~\ref{fig:honey_kin_snaps} shows a snapshot of the kinematics at $\bar{u}_2 = 1$~mm, demonstrating very close agreement between the discrete (direct) numerical simulation at $\epsilon = 1/90$ and the homogenized model. The continuum mesh contains 3749 nodes and 3626 bilinear quads. At each load step, snapshots of the kinematic fields are recorded in the form of nodal displacements and micro-deformations $(\discrete{\mathbf{u}}, \discrete{\bm{\rchi}})$, together with the corresponding reaction forces and reaction couples on the top side. Here, the nodal micro-deformations $\discrete{\bm{\rchi}}$ are given in terms of discrete out-of-plane micro-rotations $\discrete{\bm{\theta}}_3$, $(\mathrm{skw}\,\bm\rchi)_{\alpha\beta}=e_{\beta\alpha}\theta_3$. The generalized strains $(\eps_e, \gam_e)$ then follow locally from equations~\eqref{eq:epsilon_disc}--\eqref{eq:gamma_disc}.

\begin{figure}[t!]
    \centering
    \begin{subfigure}[h]{0.45\textwidth}
        \includegraphics[width=\textwidth]{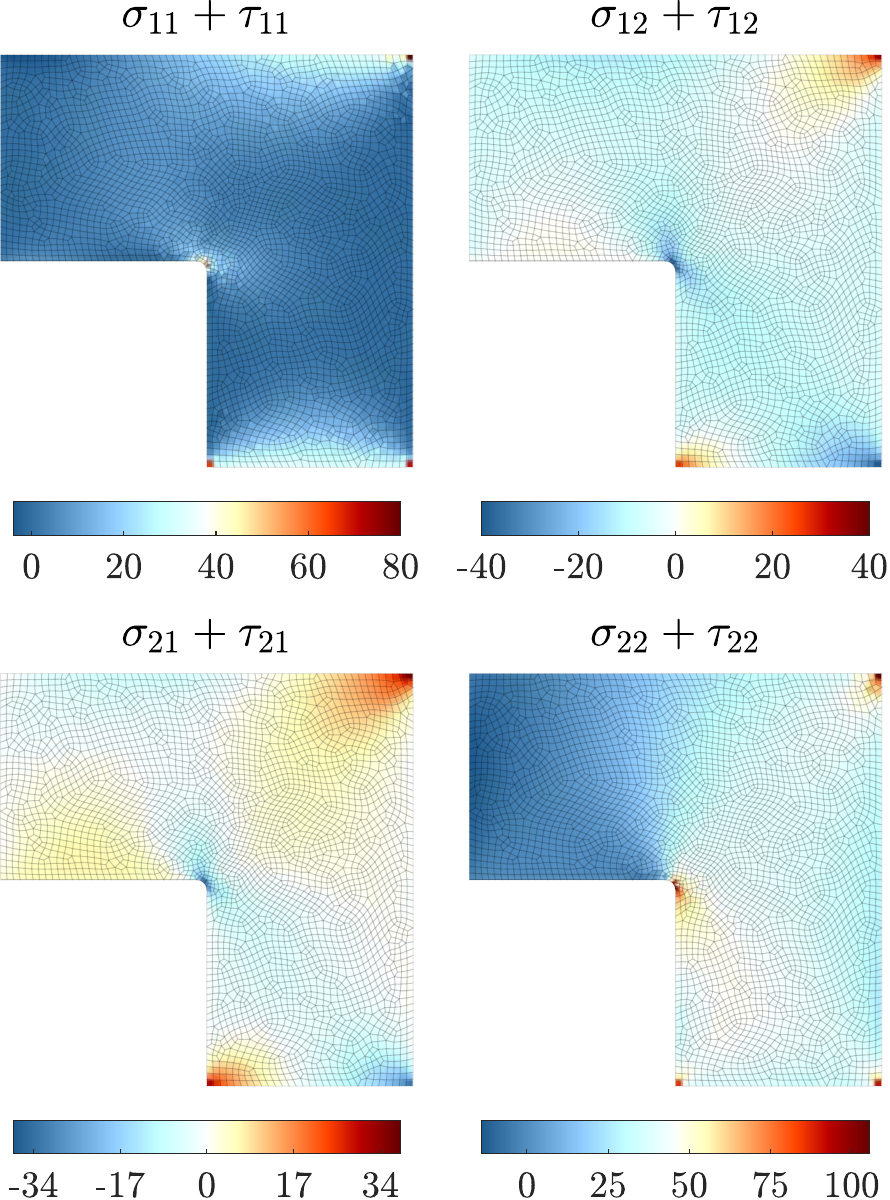}
        \subcaption{Reference}
    \end{subfigure}
    \hspace{0.5em}
    \vline
    \hspace{1em}
    \begin{subfigure}[h]{0.45\textwidth}
        \includegraphics[width=\textwidth]{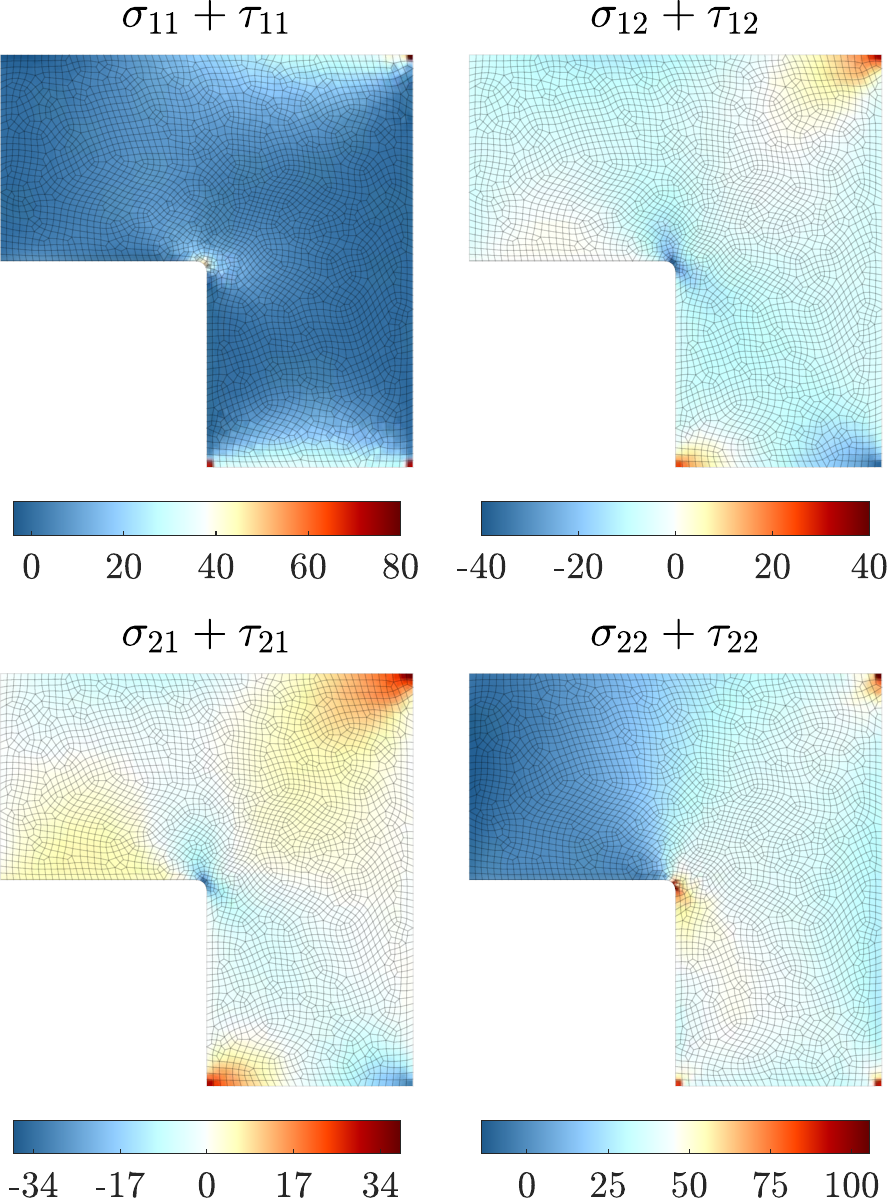}
        \subcaption{Identified}
    \end{subfigure}
    \caption{Comparison of generalized second-order stresses $\sigma_{ij}+\tau_{ij}$ [kPa] in the honeycomb L-shaped specimen, showing (a) the homogenized model simulation and (b) the identified stress fields.}
    \label{fig:honey_s_snaps}
\end{figure}

We pass the recorded data to the identification solver to obtain the corresponding stress fields and material dataset. \tB{The kinematic data is obtained either from the homogenized model simulation or by interpolating the discrete lattice displacements and rotations to the finite element mesh. The metric tensors are chosen according to~\eqref{eq:mettens}, with numerical (hyper)parameters different from those of the continuum limit: $\lambda = 50$~MPa, $\mu = 1$~MPa, and $c = 1$. Note that the metric tensor~\eqref{eq:mettens_2} includes a symmetric contribution, as we need not assume purely micropolar behavior a priori. However, because the recorded kinematic data is micropolar, the symmetric part of the relative strains $\gam_e$ vanishes; thus, only the skew-symmetric relative stress response is rendered relevant.} Again, recall that optimal metric parameters may also be chosen objectively (\ref{sec:metric}). Finally, we set the number of data points to be identified as $\bar{N} = 4352$, with a ratio of mechanical--material states $(M N_\mathrm{L})/\bar{N} \approx  100$.

\begin{figure}[t!]
    \centering
        \includegraphics[width=0.75\textwidth]{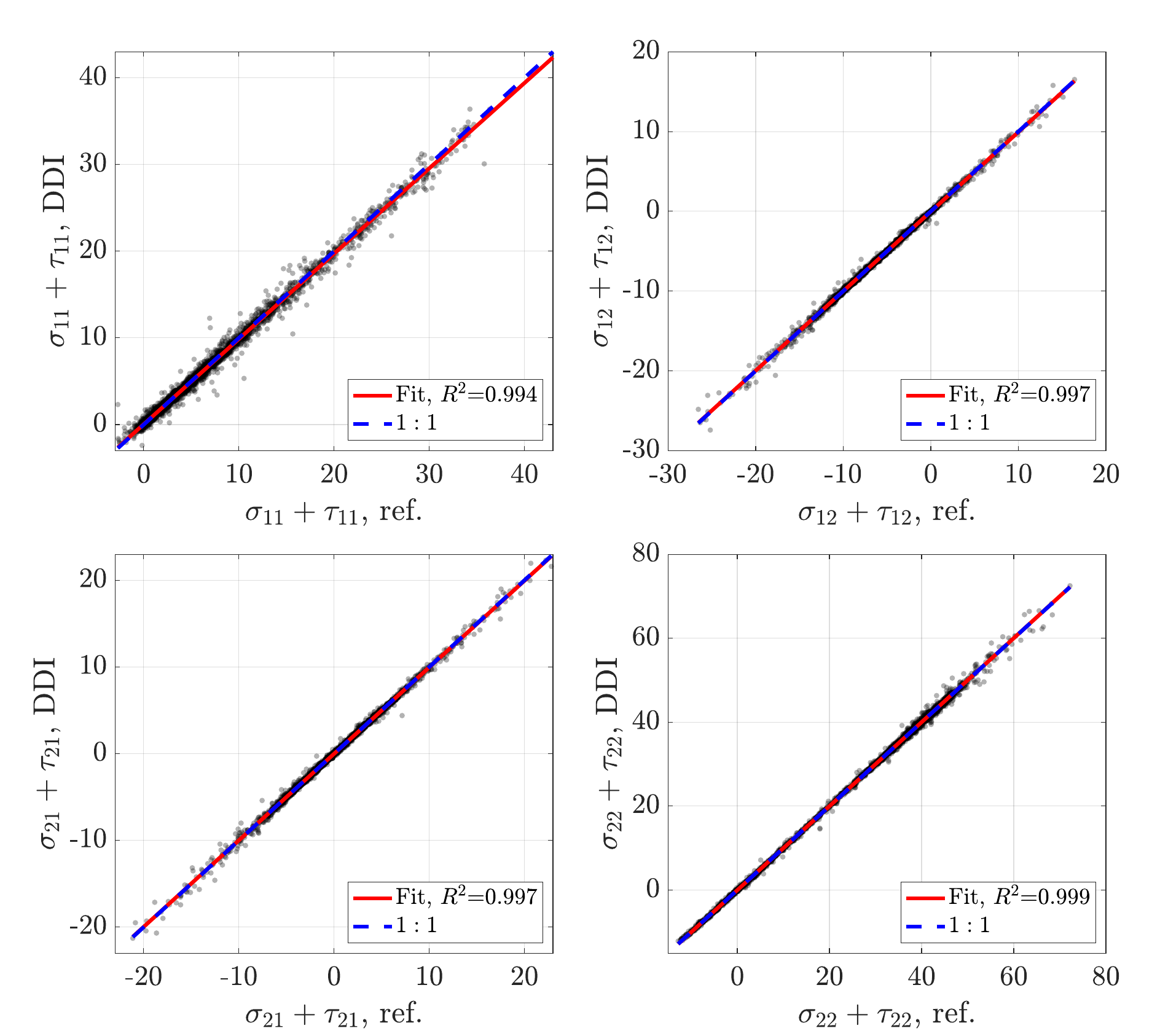}
    \caption{\tB{Correlation between identified generalized stress components and the reference homogenized model simulation.}}
    \label{fig:corr_cont_ref}
\end{figure}

\begin{figure}[t!]
    \centering
        \includegraphics[width=0.75\textwidth]{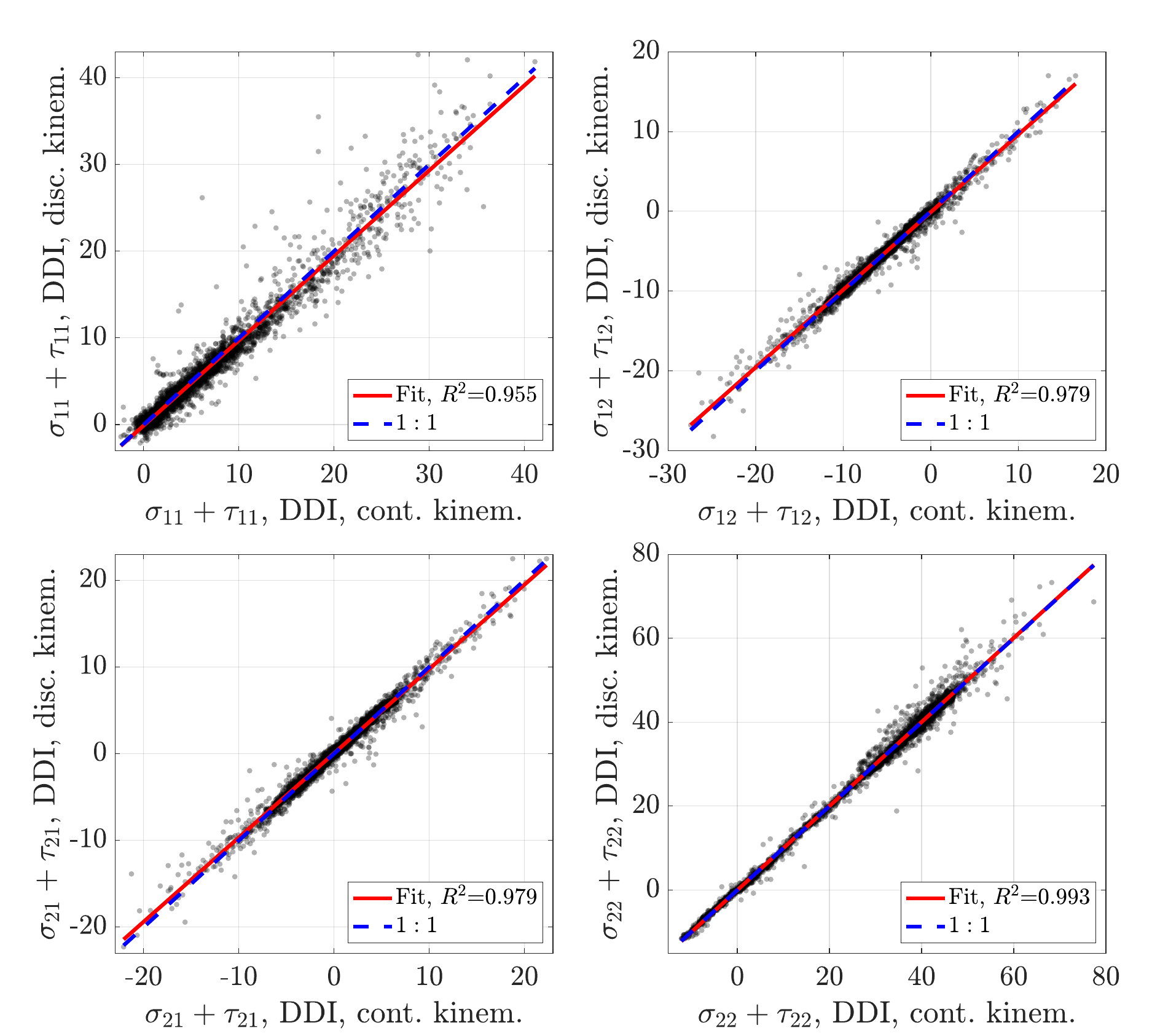}
    \caption{\tB{Correlation between identified generalized stress components obtained via discrete and continuum kinematic~data.}}
    \label{fig:corr_cont_disc}
\end{figure}

Figure~\ref{fig:honey_s_snaps} shows the generalized stress fields $\sigma_{ij}+\tau_{ij}$ at $\bar{u}_2 = 1$~mm, comparing the identified mechanical states to the reference homogenized model. We observe excellent agreement throughout the domain for all components, accurately capturing the asymmetry of $\sigma_{ij}+\tau_{ij}$ induced by the body couple. \tB{This agreement is further quantified in Figure~\ref{fig:corr_cont_ref}, where the identified generalized stress components exhibit an almost one-to-one correspondence with the reference solution. Moreover, using the kinematic data from either the homogenized model simulation or the discrete lattice leads to identified stress fields that are visually indistinguishable. Nevertheless, the latter case exhibits a modest increase in scatter (Figure~\ref{fig:corr_cont_disc}), attributed to noise amplification in the strain measures computed from the discrete kinematics. Applying regularization methods for the input gradients and preprocessing strategies for incomplete or imperfect input data~\cite{dalemat2023robustness} may help alleviate these effects in future~work.}

\begin{table}[t!]
    \caption{Reference linear relations and summary statistics of the identified data and data-driven predictions: $R^2$ values from least-squares regression and relative slope errors with respect to the homogenized parameters. The spread in the polar ratio $(\tau_{21} - \tau_{12})/(\gamma_{21} - \gamma_{12})$ is quantified using the normalized median absolute deviation (NMAD), relative to the reference~$\kappa$.}
    \centering
    \renewcommand{\arraystretch}{1.1}
    \begin{tabular}{llcccc}
    \toprule
    \multirow{2}{*}{\centering Relation} & \multirow{2}{*}{\centering Ref.~slope} &
    \multicolumn{2}{c}{Identification (DDI)} &
    \multicolumn{2}{c}{Prediction (DDCM)} \\
    & & Fit & Error [\%] & Fit & Error [\%] \\
    \cmidrule(lr){1-2}
    \cmidrule(lr){3-4}
    \cmidrule(lr){5-6}
    $\vphantom{\dfrac{\sum}{\sum}}(\sigma_{11}+\sigma_{22})/(\varepsilon_{11}+\varepsilon_{22})$
    & $2(\lambda+\mu)$
    & $R^2 = 0.99$ & 2.90 & $R^2 = 0.88$ & 6.55 \\[4pt]
    $\sigma_{12}/\varepsilon_{12}$
    & $2\mu$
    & $R^2 = 0.99 $ & 0.73 & $R^2 = 0.98$ & 0.99 \\[4pt]
    $(\tau_{21}-\tau_{12})/(\gamma_{21}-\gamma_{12})$
    & $\kappa$
    & $\mathrm{NMAD}=0.76\%$ & 0.01 & $\mathrm{NMAD}=  2.05 \%$ & 0.21 \\[4pt]
    \bottomrule
    \end{tabular}
    \label{tab:honeyDDI_all_scatter}
\end{table}

\begin{figure}[b!]
    \centering
        \includegraphics[width=0.95\textwidth]{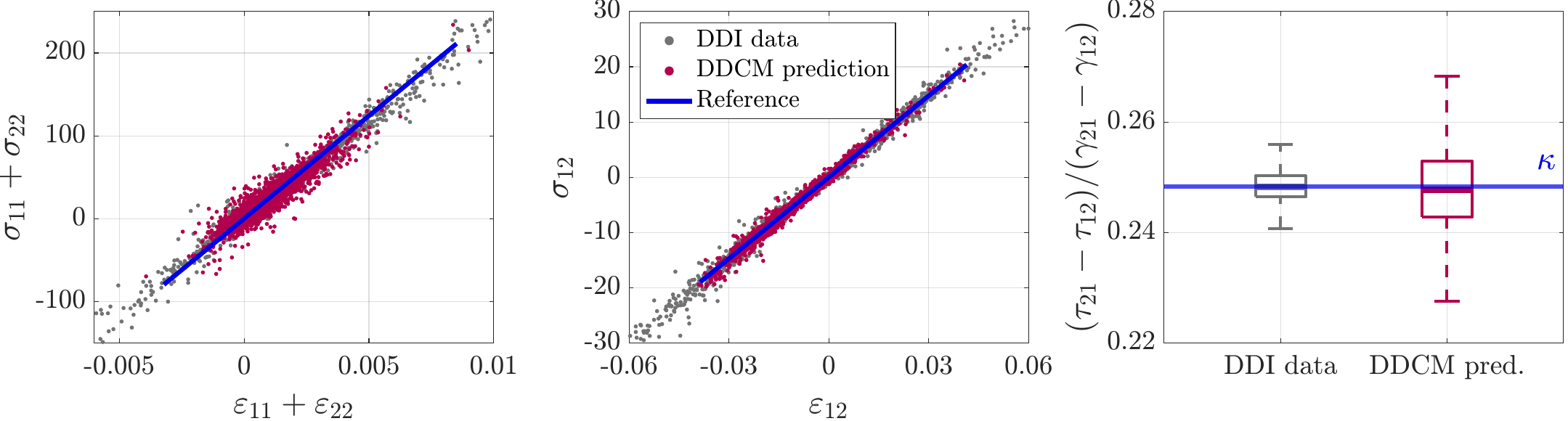}
    \caption{Material data identified from all snapshots of the honeycomb L-shaped specimen (gray, DDI data), data-driven predictions in the honeycomb double-notched specimen (red, DDCM prediction), and expected linear behavior from the homogenized model (blue). Stress quantities are shown in kPa. See tables~\ref{tab:honeyDDI_all_scatter} and~\ref{tab:honeyDDI_dists} for representative statistics.}
    \label{fig:honey_All_scatter}
\end{figure}

The identified material data points $\{\bar{\mathbf{z}}_i\}$ are shown together with the expected linear relations in Figure~\ref{fig:honey_All_scatter}. The identified data (gray dots) aligns closely with the homogenized model. The spread around the reference, most pronounced for the $\sigma_{11}+\sigma_{22}$ component, reflects the heterogeneity of the input fields. Table~\ref{tab:honeyDDI_all_scatter} reports the correlation values for the linear pairs and the relative errors in the corresponding relations. As expected from the stress contour plots, high $R^2$ values are observed for $\sigma_{11}+\sigma_{22}$ and $\sigma_{12}$, properly capturing the linear response. Moreover, the identified ratios $(\tau_{21}-\tau_{12})/(\gamma_{21}-\gamma_{12})$ are tightly distributed around the reference polar modulus $\kappa$, demonstrating accurate recovery of the micropolar constitutive response. We thus conclude that the identified relations arising from least-squares fitting of the data closely approximate all homogenized material properties.

\subsubsection{Data-driven prediction}

The identified material data from the honeycomb L-shaped specimen is now used to predict the behavior of a double-notched specimen made of the same material (Figure~\ref{fig:bvp2}b). For this purpose, we employ the micromorphic data-driven simulation framework~\cite{ulloa2024} (section~\ref{sec:fem}), providing an assessment of the reliability of the identified material data for enabling model-free predictions in a generalized continuum.

\begin{figure}[t!]
    \vspace{1em}
    \centering
    \begin{subfigure}[t]{0.45\textwidth}
        \includegraphics[width=\textwidth]{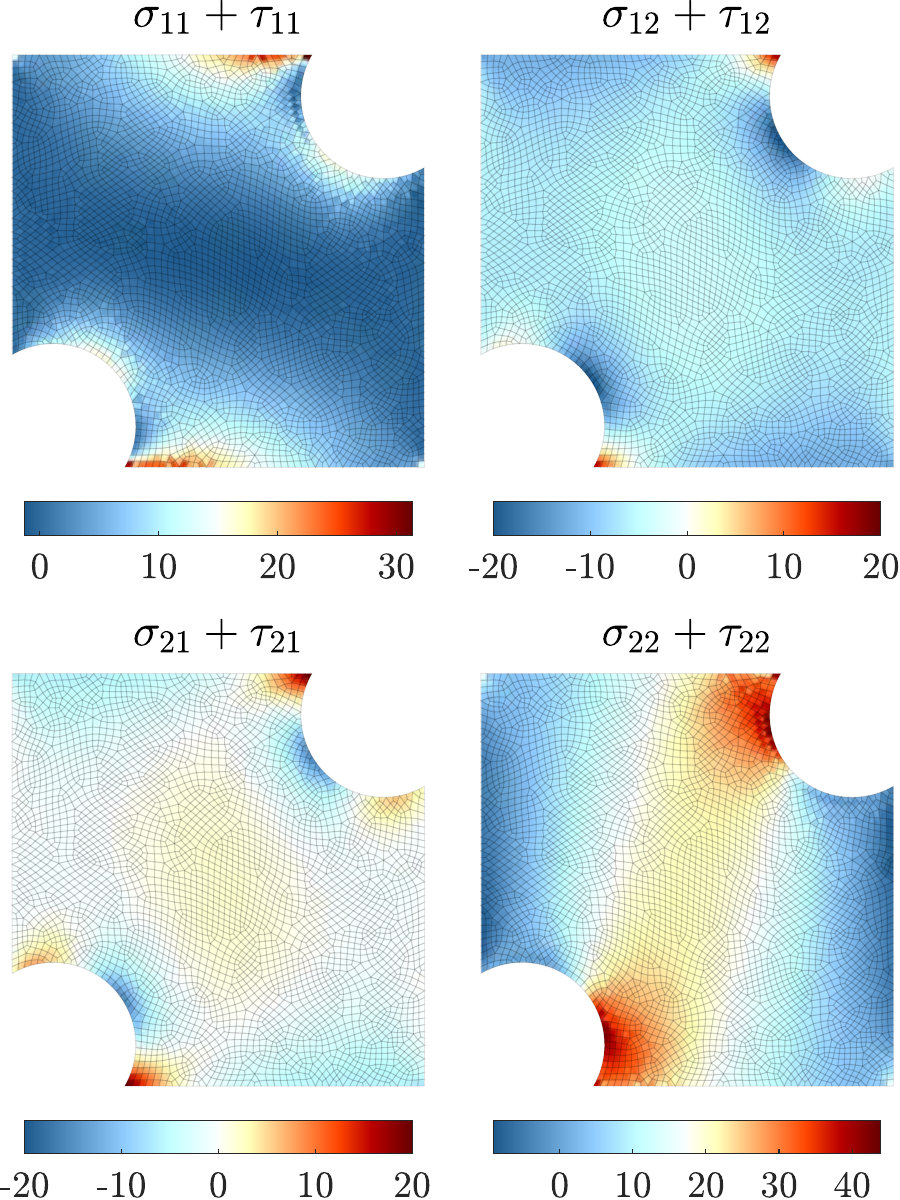}
        \subcaption{Reference}
    \end{subfigure}
    \hspace{0.5em}
    \vline
    \hspace{1em}
    \begin{subfigure}[t]{0.45\textwidth}
        \includegraphics[width=\textwidth]{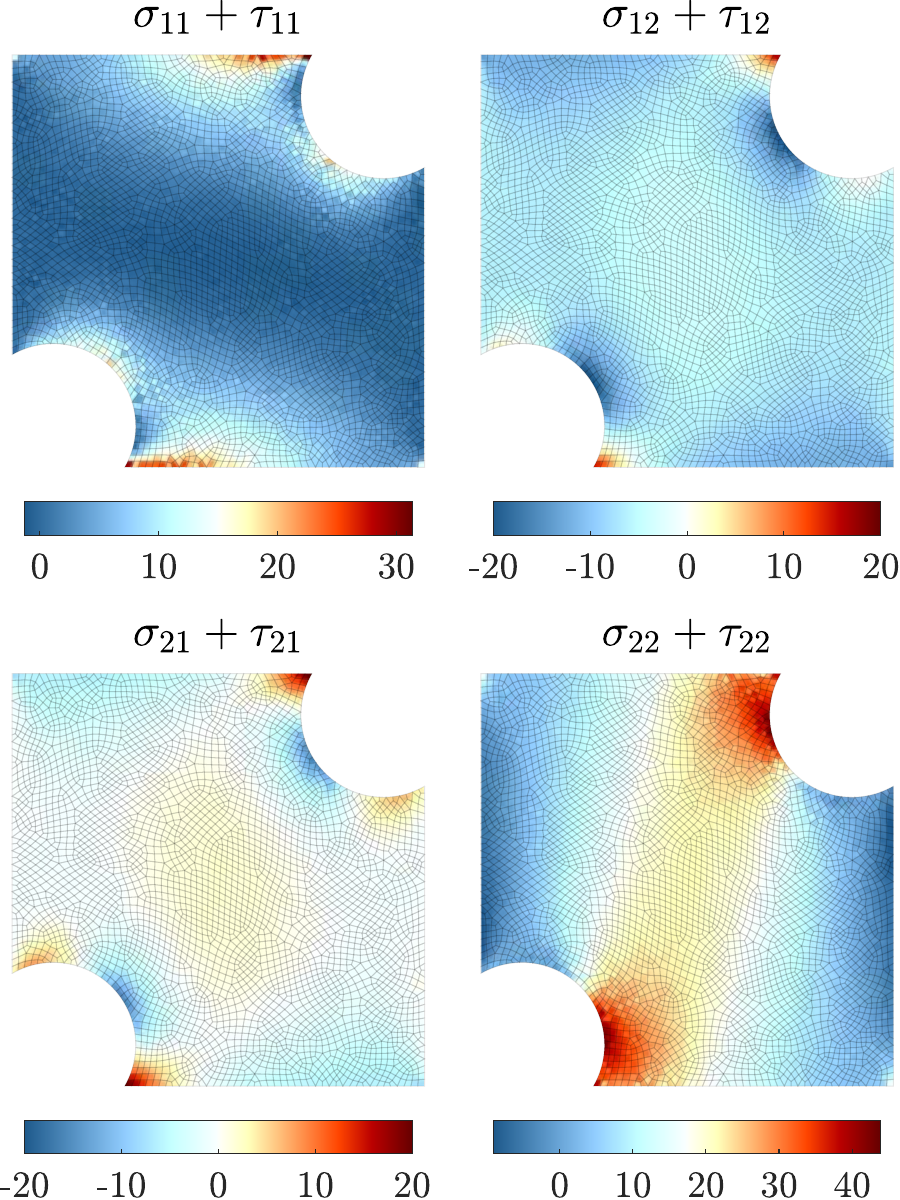}
        \subcaption{Predicted}
    \end{subfigure}
    \caption{Comparison of generalized stresses $\sigma_{ij}+\tau_{ij}$ [kPa] in the honeycomb double-notched specimen, showing (a) the homogenized model simulation and (b) the \emph{predicted} mechanical stress fields.}
    \label{fig:honey_s_pred}
\end{figure}

\begin{figure}[t!]
    \centering
    \includegraphics[width=0.8\textwidth]{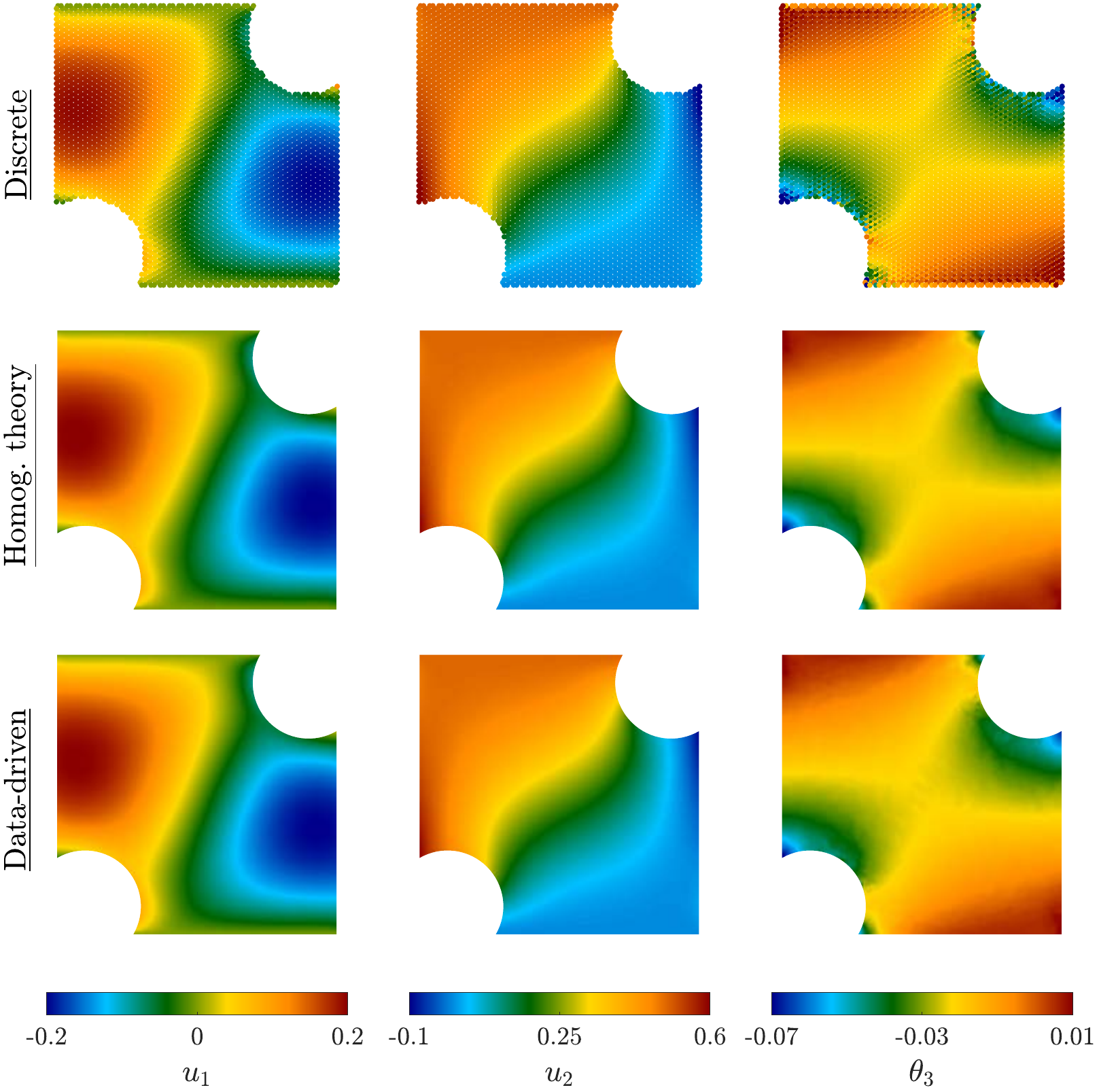}
    \caption{Kinematic fields in the honeycomb double-notched specimen: displacements $u_i$ [mm] and micro-rotations $\theta_3$, showing the direct numerical simulation of the discrete metastructure (top), the homogenized model (middle), and the DDI-informed data-driven predictions (bottom).}
    \label{fig:honey_kin_pred}
\end{figure}

Figure~\ref{fig:honey_s_pred} shows the generalized stress fields $\sigma_{ij}+\tau_{ij}$ at $\bar{u}_2 = 0.35$~mm, comparing the \emph{predicted} mechanical states to the reference homogenized model in the double-notched specimen. As for the identified states in the L-shaped specimen, we observe an excellent agreement in all four components. This behavior is also observed in Figure~\ref{fig:honey_All_scatter}, where the predicted generalized stress--strain states lie, as expected, within the DDI data points and are closely distributed around the reference response. The tendency toward the expected linear relations is confirmed by the statistics in Table~\ref{tab:honeyDDI_all_scatter}. Furthermore, Figure~\ref{fig:honey_kin_pred} presents the kinematic fields $(u_1,u_2,\theta_3)$, comparing the direct numerical simulation of the discrete metastructure to the homogenized model and the data-driven predictions. Again, we observe a notable agreement between the DDI-informed data-driven simulation and the expected response.

\begin{table}[t!]
    \caption{Metric distances of the phase-space coordinates $\mathbf{z} = (\eps, \gam, \sig, \btau)$ in the honeycomb double-notched specimen, reporting (i) the relative mechanical--material errors in the data-driven simulation (DDCM) and (ii) the relative errors with respect to the reference homogenized model.}
    \centering
    \renewcommand{\arraystretch}{1.3}
    \begin{tabular}{llcc}
    \toprule
    Variable & Metric $d(\cdot)$ 
    & DDCM, $100\dfrac{d(\Box - \bar{\Box})}{d(\Box)}$ [\%]
    & Reference, $100\dfrac{d(\Box - \Box_{\mathrm{ref}})}{d(\Box_{\mathrm{ref}})}$ [\%]\\
    \midrule
    Strain $\eps$      
    & $\Vert\cdot\Vert_{\Ctens}$ 
    & 6.35 & 4.10 \\
    Stress $\sig$      
    & $\Vert\cdot\Vert_{\Ctens^{-1}}$ 
    & 5.88 & 6.81 \\
    Rel.\ strain $\gam$ 
    & $\Vert\cdot\Vert_{\Dtens}$ 
    & 4.06 & 1.81 \\
    Rel.\ stress $\btau$ 
    & $\Vert\cdot\Vert_{\Dtens^{-1}}$
    & 0.69 & 3.76 \\
    All $\mathbf{z}$    
    & $\Vert\cdot\Vert$ 
    & 5.61 & 5.19 \\
    \bottomrule
    \end{tabular}
    \label{tab:honeyDDI_dists}
\end{table}

The predictive accuracy is further assessed quantitatively in Table~\ref{tab:honeyDDI_dists}. Specifically, we report the global distance metric~\eqref{eq:globmet} and its individual contributions from each phase-space variable, namely
$$
    \Vert \mathbf{z} \Vert^2 = \underbrace{\frac{1}{2} \sum_{e=1}^M w_e\,\Ctens_e\,\eps_e\cdot\eps_e}_{\displaystyle\Vert \eps \Vert_{\Ctens}^2} + \underbrace{\frac{1}{2} \sum_{e=1}^M w_e\,\Ctens^{-1}_e\,\sig_e\cdot\sig_e}_{\displaystyle\Vert \sig \Vert_{\Ctens^{-1}}^2} + \underbrace{\frac{1}{2} \sum_{e=1}^M w_e\,\Dtens_e\,\gam_e\cdot\gam_e}_{\displaystyle\Vert \gam \Vert_{\Dtens}^2}  + \underbrace{\frac{1}{2} \sum_{e=1}^M w_e\,\Dtens^{-1}_e\,\btau_e\cdot\btau_e}_{\displaystyle\Vert \btau \Vert_{\Dtens^{-1}}^2},
$$
considering (i) the errors between the predicted material solution $\bar{\mathbf{z}}$ and the mechanical solution $\mathbf{z}$, and (ii) the errors between the predicted mechanical solution $\mathbf{z}$ and the reference solution $\mathbf{z}_\mathrm{ref}$ from homogenization theory. All error metrics remain well below 10\%, confirming the reliability of the identified data for predicting the response of unseen BVPs in a generalized continuum.

\section{Summary and concluding remarks}
\label{sec:conclusion}

We have presented a data-driven identification (DDI) framework to extract material data from full-field kinematics and force measurements in generalized continua at small strains. \tB{The central question addressed is whether this input information, presumed observable, can reveal the active generalized stress--strain components and their constitutive response within a prescribed micromorphic phase space, without prescribing closure relations.} By adopting a model-free data-driven identification approach, we have relied on (i) non-classical compatibility and balance equations, and (ii) the availability of full-field kinematic data\md displacements and micro-deformations\md and the corresponding applied forces to infer the associated generalized stress fields along with a representative material dataset. Notably, our approach circumvents the definition of heuristic kinematics and boundary conditions typically required in RVE-based micromorphic homogenization~methods. The identified data is suitable for either calibration of constitutive models or model-free data-driven simulations of generalized continua. 

\tB{The reliability of the identification process was first validated via synthetic studies, considering both a full elastic micromorphic material and an elastoplastic microstrain material. These examples demonstrated the recovery of non-trivial features, such as non-symmetric and higher-order stress components, as well as both linear and nonlinear constitutive responses for active phase-space components.} We then presented an application to mechanical metamaterials, including predictions in unseen BVPs of a honeycomb metastructure. In this case, the ground-truth behavior used for benchmarking is known explicitly from homogenization theory and corresponds to a zeroth-order micropolar continuum. \tB{The identified material data showed excellent agreement with this reference. Moreover, the identification framework remains applicable beyond the tested lattice geometry, in cases where closed-form homogenized coefficients may be unavailable but discrete kinematics and applied forces can still be~obtained.}

\tB{Future work should address systematic selection of metric parameters and dataset size, regularization of noisy or discrete kinematic fields, and extensions to material instabilities and finite deformations. We further emphasize that the present approach remains agnostic with respect to constitutive closure relations, while the micromorphic phase space is assumed as a point of departure.} This approach is sufficiently general to encompass a wide range of continuum behaviors as special cases, including continuum limits with provable convergence from discrete systems (e.g., the example in section~\ref{sec:MMM}). Ultimately, the specific form of the continuum should be established through such rigorous homogenization analyses, which may reveal the relevant phase-space quantities without necessarily yielding explicit, closed-form parameters. In those scenarios, the present data-driven framework provides a flexible means of identifying generalized material data that may otherwise not be directly accessible.

\section*{Acknowledgements}

JU gratefully acknowledges support from the Digital Twins seed funding initiative provided by the Department of Mechanical Engineering at the University of Michigan.
LS gratefully acknowledges the financial support of Nantes Université excellence program NExT, through the funding of International Research Partnership project iDDrEAM.

\appendix

\section{Metric tensors}
\label{sec:metric}

Let us recall that, in the local distance metric~\eqref{eq:locmet}, the phase-space coordinates
$$(\eps_e,\sig_e)\in\mathbb{R}^{n_\varepsilon} \times\mathbb{R}^{n_\varepsilon}, \qquad (\gam_e,\btau_e)\in\mathbb{R}^{n_\chi}\times\mathbb{R}^{n_\chi}, \qquad (\zet_e,\bmu_e)\in\mathbb{R}^{n_\zeta}\times\mathbb{R}^{n_\zeta}$$
are Voigt representations of the generalized stress--strain states 
$$(\eps,\sig)\in\mathbb{R}^{n\times n}_{\mathrm{sym}}\times\mathbb{R}^{n\times n}_{\mathrm{sym}}, \qquad (\gam,\btau)\in\mathbb{R}^{n\times n}\times\mathbb{R}^{n\times n}, \qquad (\zet,\bmu)\in\mathbb{R}^{n\times n\times n}\times\mathbb{R}^{n\times n\times n}.$$ 
Hence, the metric operators  $\Ctens_e\in\mathbb{R}^{n_\varepsilon\times n_\varepsilon}_{\mathrm{sym},+}$, $\Dtens_e\in\mathbb{R}^{n_\chi\times n_\chi}_{\mathrm{sym},+}$, and $\Atens_e\in\mathbb{R}^{n_\zeta\times n_\zeta}_{\mathrm{sym},+}$ also correspond to Voigt representations of tensors $\Ctens$ (fourth-order), $\Dtens$ (fourth-order), and $\Atens$ (sixth-order). 

In standard isotropic constitutive models, $\Ctens$ owns the usual minor and major symmetries $\mathsf{C}_{ijkl}=\mathsf{C}_{jilk}=\mathsf{C}_{klij}$, while, owing to the non-symmetry of $\bm{\rchi}$, in general, $\Dtens$ and $\Atens$ only possess major symmetries $\mathsf{D}_{ijkl}=\mathsf{D}_{klij}$ and $\mathsf{A}_{ijklmn}=\mathsf{A}_{lmnijk}$. However, the data-driven approach treats these as numerical operators. As such, we choose simple forms with the following structure:
\begin{subequations}\label{eq:mettens}
\begin{align}
     &\mathsf{C}_{ijkl} \;=\; \lambda \, \delta_{ij}\delta_{kl} \;+\; \mu\, (\delta_{ik}\delta_{jl} + \delta_{il}\delta_{jk}),   
      \label{eq:mettens_1}\\
     &\mathsf{D}_{ijkl} \;=\; c\,\mathsf{C}_{ijkl} \;+\; c\,\mu\, (\delta_{ik}\delta_{jl} - \delta_{il}\delta_{jk}),
     \label{eq:mettens_2}\\
     &\mathsf{A}_{ijklmn} \;=\; \ell^2\,\mathsf{C}_{ijlm}\,\delta_{kn} \;+\; \mu\,\ell^2 \, (\delta_{il}\delta_{jm} - \delta_{im}\delta_{jl})\,\delta_{kn}. \label{eq:mettens_3}
\end{align}
\end{subequations}
The generalized phase space is then metrized by \tB{four positive metric constants, $\bm{\theta}\coloneqq(\lambda,\mu,c,\ell)$}. 

A convenient strategy to define these coefficients objectively, adopted in~\citet{karapiperis2021b} for forward data-driven problems of micropolar continua, is to minimize the distance~\eqref{eq:locmet}, locally or globally, at fixed mechanical and material states $(\mathbf{z}_e,\bar{\mathbf{z}}_e)$ but varying metric constants. For instance, for the global case,
\begin{equation}
   \bm{\theta} \in \argmin_{\displaystyle\tB{\tilde{\bm{\theta}}\in\mathbb{R}^4_+}}\,  \Vert \bar{\mathbf{z}}-\mathbf{z}\Vert^2\,(\tilde{\bm{\theta}}),
\end{equation}
where we write the dependence of the global distance metric on the metric constants explicitly. This procedure is not employed here for the sake of simplicity but can be incorporated in Algorithm~\ref{alg:fixedpoint} to select optimal metric tensors~\eqref{eq:mettens} objectively.

\small
\bibliography{literature}

\end{document}